\documentstyle[12pt]{article}
\def\be{\begin{equation}}
\def\ee{\end{equation}}

\def\hang{\hangindent\parindent}
\def\textindent#1{\indent\llap{#1\enspace}\ignorespaces}
\def\re{\par\hang\textindent}

\begin{document}

\vspace*{1cm}

\centerline{\large\bf Zeros of the Riemann zeta function on the
critical line} \vskip 0.6cm \centerline{Shaoji Feng}
\centerline{Academy of Mathematics and Systems Science}
\centerline{Chinese Academy of Sciences} \centerline{Beijing 100190}
\centerline{P. R. China} \centerline{email: fsj@amss.ac.cn}

\vskip 1cm

{\small \noindent{\bf Abstract.} We introduce a new mollifier and
apply the method of Levinson and Conrey to prove that at least
$41.28\%$ of the zeros of the Riemann zeta function are on the
critical line. The method may also be used to improve other results
on zeros relate to the Riemann zeta function, as well as conditional
results on prime gaps.}

\noindent {\bf 2000 Mathematics Subject Classification.} 11M26,
11M06

\noindent {\bf Keywords.} Riemann zeta function, zeros, critical
line, mollifier

\vskip 8mm

\def\theequation{1.\arabic{equation}}
\setcounter{equation}{0} \centerline{\large\bf 1. \ Introduction}
\vskip 5mm

\noindent  Let $\zeta(s)$ denote the Riemann zeta function, where
$s=\sigma+it$. It is defined for $\sigma>1$ by
$$ \zeta(s)=\sum_{n=1}^{\infty}n^{-s}=\prod_p(1-p^{-s})^{-1},$$
where $p$ runs over the prime numbers, and has a meromorphic
continuation to the whole complex plane with its only pole, a simple
pole at $s=1$. It satisfies the functional equation \be
\xi(s)=\xi(1-s),\ee where the entire function $\xi(s)$ is defined by
\be \xi(s)=H(s)\zeta(s)\ee with \be
H(s)=\frac{s(s-1)}{2}\pi^{-s/2}\Gamma(\frac{s}{2}).\ee

Let $N(T)$ denote the number of zeros of $\zeta(s)$, $s=\sigma+it$,
in the rectangle $0<\sigma<1,\ 0<t\leq T$, each zeros is counted
with multiplicity, von Mangoldt proved that (see [24])
$$ N(T)=\frac{T}{2\pi}(\log\frac{T}{2\pi}-1)+\frac78+S(T)+O(\frac{1}{T}),$$
$$S(T)=\frac{1}{\pi}\arg\zeta(\frac12+iT)=O(\log T), \ \mbox{as}\ \
T\rightarrow\infty.$$ Let $N_0(T)$ be the number of zeros of
$\zeta(\frac12+it)$ on $0<t\leq T$, each zeros is counted with
multiplicity, $N_{0s}(T)$ be the number of simple zeros of
$\zeta(\frac12+it)$ on $0<t\leq T$. The Riemann Hypothesis says that
$N_0(T)=N(T)$, the Simple Zeros Conjecture combined with the Riemann
Hypothesis says that $N_{0s}(T)=N_0(T)=N(T)$.

It was proved for the first time by Hardy [11] in 1914 that
$\zeta(s)$ has infinitely many zeros on the critical line
$\sigma=\frac12$, thus $$ N_0(T)\rightarrow\infty\ \ \ \mbox{as}\ \
\ T\rightarrow\infty.$$ Hardy's qualitative result was given a
quantitative form
$$ N_0(T)\geq AT $$ for some $A>0$ and $T$ large enough, by Hardy
and Littlewood [12] in 1921, and later, with an explicit value of
$A$, the same result was obtaned by Siegel [23] in 1932 with a
rather different method.

Let \be \kappa=\liminf_{T\rightarrow\infty}\frac{N_0(T)}{N(T)},\ee
\be \kappa^*=\liminf_{T\rightarrow\infty}\frac{N_{0s}(T)}{N(T)}.\ee
In 1942, Selberg [21] proved that there is an effectively computable
positive constant $A$ such that
$$ \kappa\geq A.$$  Selberg's proof involved combing a `mollifier'
to compensate for irregularities in the size of $|\zeta(s)|$ and the
method of Hardy and Littlewood.

In 1974, Levinson [14] combined Siegel's idea and Selberg's idea and
proved that
$$ \kappa\geq 0.3420.$$

The Levinson method involve the following main issues:

i) `Perturb' the Riemann zeta function $\zeta(s)$ to a function
$f(s)$ such that (if the Riemann Hypothesis is true) $f(s)$ remains
has no (or few) zeros in the rectangle $\frac12<\sigma<1,\ 0<t\leq
T$ and the zeros of $\zeta(s)$ on the critical line are `pushed' to
the left of the critical line. Levinson's choice of $f(s)$ in [14]
is essentially $\zeta(s)+\frac{\zeta'(s)}{\log T}$.

ii) Use the Littlewood Lemma and a convexity inequality to estimate
the upper bound of zeros $\beta+i\gamma$ of $f(s)$ in the rectangle
$\frac12-\frac{R}{\log T}<\sigma<1,\ 0<t\leq T$ with weight
$\beta-\frac12+\frac{R}{\log T}$, and then get an upper bound
estimation of zeros of $f(s)$ in the rectangle $\frac12<\sigma<1,\
0<t\leq T$.

iii) A mollifier $\psi(s)$ is used in ii) to compensate for
irregularities in the size of $|f(s)|$. The mollifier Levinson
constructed in [14] is \be \psi(s)=\sum_{1\leq j\leq
y}\frac{\mu(j)}{j^{\frac{R}{\log T}}j^s}\frac{\log y/j}{\log y},\ee
where $y=T^{\theta}$ (with $\theta=\frac12-\varepsilon$) is the
length of the mollifier.

iv) Since (see Levinson and Montgamery [17]) unconditionally
$\zeta(s)$ and $f(s)$ have almost the same number of zeros in the
rectangle $\frac12<\sigma<1,\ 0<t\leq T$, ii) gives an upper bound
estimation of zeros of $\zeta(s)$ in the rectangle
$\frac12<\sigma<1,\ 0<t\leq T$ and then a lower bound estimation of
$N_0(T)$.

In an attempt to semi-optimize the coefficients of the mollifier
$\psi(s)$, Levinson [15] choose  \be \psi(s)=\sum_{1\leq j\leq
y}\frac{\mu(j)}{j^{\frac{2R}{\log T}}j^s}\frac{y^{\frac{2R}{\log
T}}-j^{\frac{2R}{\log T}}}{y^{\frac{2R}{\log T}}-1}\ee with
$y=T^{\frac12-\varepsilon}$, and gives the result $$ \kappa\geq
0.3474.$$

Heath-Brown [13] and Selberg independently noticed that the zeros
located by Levinson's method are simple zeros of $\zeta(s)$, thus
$$ \kappa^*\geq
0.3474.$$ It is the first time one proved unconditionally that there
are infinitely many simple zeros of $\zeta(s)$ in the critical
strip.

Lou and Yao [18] choose the mollifier \be \psi(s)=\sum_{1\leq j\leq
y}\frac{\mu(j)}{j^{\frac{2R}{\log T}}j^s}\frac{y^{\frac{2R}{\log
T}}-j^{\frac{2R}{\log T}}}{y^{\frac{2R}{\log T}}-1}+h_0\sum_{1\leq
j\leq y}\frac{\mu(j)}{j^{\frac{R}{\log T}}j^s}\frac{\log \frac{y}{j}
\log j} {\log^2y}\ee with $h_0$ a real constant and
$y=T^{\frac12-\varepsilon}$. This leads to
$$ \kappa\geq 0.3484.$$ Lou also announced a result of $ \kappa\geq
0.35$, but without detailed proof.

In [4], Conrey use the more general mollifier \be
\psi(s)=\sum_{1\leq j\leq y}\frac{\mu(j)}{j^{\frac{R}{\log
T}}j^s}P(\frac{\log y/j}{\log y}),\ee where
$y=T^{\frac12-\varepsilon}$ and $P$ is an analytic function with
$P(0)=0$ and $P(1)=1$. The function $f(s)$ is also generalized to
$$ f(s)=Q(-\frac{1}{\log T}\frac{d}{ds})\zeta(s),$$
where $Q$ is a real polynomial with $Q(0)=1$ and $Q'(x)=Q'(1-x)$. By
choosing $P$ and $Q$ appropriately this gives
$$ \kappa\geq 0.3658.$$

With Levinson's original $f(s)$ and the mollifier (1.9), Conrey [5]
obtain
$$ \kappa^*\geq 0.3485.$$

Anderson [1] use the mollifier (1.7) and $$
f(s)=\zeta(s)+a_1\frac{\zeta'(s)}{\log T},$$ where $a_1$ is an
arbitrary real number. This gives
$$ \kappa^*\geq 0.3532.$$

In 1989, Conrey [6] proved that
$$ \kappa\geq 0.4088$$
and
$$ \kappa^*\geq 0.4013.$$ These significant improvement has been
obtained by using a mollifier of length $y=T^{\theta}$ with
$\theta=\frac47-\varepsilon$. The work of Deshouillers and Iwaniec
[7,8] on Kloosterman sums is used to estimate the error terms of the
mean value integral for this longer mollifier.

Recently, in [3], Bui, Conrey and Young introduced a two-piece
mollifier and proved
$$ \kappa\geq 0.4105,$$
$$ \kappa^*\geq 0.4058.$$

Although the Levinson method and its modification are most
successful and hopeful at present for estimation of the proportion
of zeros of the Riemann zeta function on the critical line, there
are several aspects that prevent one from getting essentially better
results:

a) Zhang [25] proved on Small Gap Zeros conjecture and the Riemann
Hypothesis and later Feng [10] proved only on Small Gap Zeros
conjecture that $\zeta'(s)$ has a positive proportion of zeros near
the critical line. By the functional equation this means that
$\zeta(s)+\frac{\zeta'(s)}{\log T}$ has a positive proportion of
zeros near the critical line. That is, although the zeros of
$\zeta(s)$ on the critical line are `pushed' to the left of the
critical line, there is a positive proportion of zeros which are not
`pushed' far away from the critical line enough. Therefore, if one
expect to prove that $100\%$ of the zeros of the Riemann zeta
function are on the critical line by using the Levinson Method, they
must let $R\rightarrow 0$ or construct an essentially different
$f(s)$.

b) Farmer [9] proposed the `$\theta=\infty$ conjecture' and proved
that this implies $100\%$ of the zeros of the Riemann zeta function
are on the critical line. This shows that the length of the
mollifier are key for the Levinson Method. However, we can see in
Conrey [6] that it is very difficult to deal with the error terms
for longer mollifier.

c) For given $f(s)$ and given length of the mollifier, it is too
complicated to optimize exactly the coefficients of the mollifier.

In this paper, we prove the following

\vskip 3mm

{\bf Theorem 1.} \be \kappa\geq 0.4128. \ee

\vskip 3mm

{\bf Remark.} One can combine our method and that of Bui, Conrey and
Young [3] to improve both $\kappa$ and $\kappa^*$.

The framework of proof follows that of Levinson and Conrey. The main
new element here is the use of a different mollifier:
\begin{eqnarray} \psi(s)&=&\sum_{j\leq
y}\frac{\mu(j)}{j^{\frac{R}{\log T}+s}}P_1(\frac{\log y/j}{\log
y})+\sum_{j\leq y_1}\frac{\mu(j)}{j^{\frac{R}{\log
T}+s}}\Big(\sum_{p_1p_2|j}\frac{\log p_1\log
p_2}{\log^2y_1}P_2(\frac{\log
y_1/j}{\log y_1})\nonumber\\
&&+\sum_{p_1p_2p_3|j}\frac{\log p_1\log p_2\log
p_3}{\log^3y_1}P_3(\frac{\log y_1/j}{\log
y_1})+\cdots\nonumber\\
&&+\sum_{p_1p_2\cdots p_I|j}\frac{\log p_1\log p_2\cdots\log
p_I}{\log^Iy_1}P_I(\frac{\log y_1/j}{\log y_1})\Big),\end{eqnarray}
where $y=T^{\theta}$ with $\theta=\frac47-\varepsilon$ and
$y_1=T^{\theta_1}$ with $\theta_1=\frac12-\varepsilon$, $I\geq 2$ is
a integer, $P_1$ is a real polynomial with $P_1(0)=0$ and
$P_1(1)=1$, $P_l (l=2,\cdots I)$ are real polynomials with
$P_l(0)=0$, $p_1,p_2,\cdots,p_I$ runs over the prime numbers.

Selberg [22] proved that for $a_1=1$, \be \sum_{m,n\leq
y}\frac{a_ma_n}{mn}(m,n)\geq 1/(\sum_{n\leq
y}\frac{\mu^2(n)}{\phi(n)})=\frac{1+o(1)}{\log y},\ee where $(m,n)$
is the greatest common divisor of $m$ and $n$, $\phi(n)$ is Euler's
phi function, and the equality are achieved when \be
a_n=\mu(n)\frac{nL(y/n;n)}{\phi(n)L(y;1)},\ee where $$
L(X,n)=\sum_{m\leq X\atop (m,n)=1}\frac{\mu^2(m)}{\phi(m)}.$$ For
fixed $n$ it is not difficult to show that $$
L(X,n)\sim\frac{\phi(n)}{n}\log X$$ and it may be expected that
$L(X,n)$ is approximately this size for most $n$. This suggests that
$a_n$ in (1.13) are approximately \be \mu(n)\frac{\log y/n}{\log
y}.\ee Thus the motivation of Levinson's choice of the coefficients
of the mollifier can be understood as to minimize the quadratic form
in (1.12), which is simplification of the main term for the mean
value integral (see Mongomery [19]).

On the other hand, the mollifier (1.6)-(1.9) can also be understood
as `continuous' truncation of the Dirichlet series of
$\frac{1}{\zeta(s)}$. Therefore the mollifier can be seemed as try
to mollify $\zeta(s)$.

However, while in Selberg's method we need to mollify $\zeta(s)$, in
Levinson's method we really need to mollify the perturbed function
$f(s)$. Our motivation of choice of (1.12) is try to mollify
$\zeta(s)+\frac{\zeta'(s)}{\log T}$. Consider the Dirichlet series
\be \frac{1}{\zeta(s)+\frac{\zeta'(s)}{\log
T}}=\frac{1}{\zeta(s)}-\frac{\zeta'(s)}{\log
T\zeta^2(s)}+\frac{\zeta'^2(s)}{\log^2
T\zeta^3(s)}-\frac{\zeta'^3(s)}{\log^3 T\zeta^4(s)}+\cdots.\ee For
$j$ be a square-free positive integer, we have
\begin{eqnarray}(\mu*\Lambda)(j)&=&-\mu(j)\log j,\nonumber\\
 (\mu*\Lambda*\Lambda)(j)&=&\mu(j)\sum_{p_1p_2|j}\log p_1\log
 p_2,\nonumber\\
(\mu*\Lambda*\Lambda*\Lambda)(j)&=&-\mu(j)\sum_{p_1p_2p_3|j}\log
p_1\log p_2\log p_3,\\
&\cdots,&\nonumber \end{eqnarray} where $f*g$ denotes the Dirichlet
convolution of arithmetic functions $f$ and $g$. For those $j$ which
contains a square divisor, the coefficients $a_j$ defined according
to (1.15) contribute a lower order term for the mean value integral.
Therefore, the mollifier (1.11) can be seemed as simplification of
`continuous' truncation of the Dirichlet series
$\frac{1}{\zeta(s)+\frac{\zeta'(s)}{\log T}}$.

We mention that the method we used here may also apply to improve
other results on zeros relate to the Riemann zeta function, as well
as conditional results on prime gaps.

 \vskip 8mm

\def\theequation{2.\arabic{equation}}
\setcounter{equation}{0} \centerline {\large\bf 2. \ Beginning of
the proof and some lemmas} \vskip 5mm

\noindent By a standard discussion as in Conrey [6] section 3 (see
Levinson [14] and Conrey [4] also), we have by (1.1), Littlewood
Lemma and the fact \be a_j\ll 1, a_1=1\ee that

\vskip 3mm

{\bf Lemma 1}. Let $T$ be a large parameter and $L=\log T$, $R$ be a
positive real number, $\sigma_0=\frac12-\frac{R}{L}$, $y=T^{\theta}$
with $\theta=\frac47-\varepsilon$ and $y_1=T^{\theta_1}$ with
$\theta_1=\frac12-\varepsilon$, \be B(s)=\sum_{j\leq
y}\frac{a_j}{j^{s+\frac{R}{L}}},\ee for $y_1<j\leq y$, \be
a_j=\mu(j)P_1(\frac{\log y/j}{\log y}),\ee for $j\leq y_1$,
\begin{eqnarray} a_j&=&\mu(j)\Big(P_1(\frac{\log y/j}{\log y})+\sum_{p_1p_2|j}\frac{\log p_1\log
p_2}{\log^2y_1}P_2(\frac{\log
y_1/j}{\log y_1})\nonumber\\
&&+\sum_{p_1p_2p_3|j}\frac{\log p_1\log p_2\log
p_3}{\log^3y_1}P_3(\frac{\log y_1/j}{\log
y_1})+\cdots\nonumber\\
&&+\sum_{p_1p_2\cdots p_I|j}\frac{\log p_1\log p_2\cdots\log
p_I}{\log^Iy_1}P_I(\frac{\log y_1/j}{\log y_1})\Big),\end{eqnarray}
where $I\geq 2$ is a integer, $P_1$ is a real polynomial with
$P_1(0)=0$ and $P_1(1)=1$, $P_l (l=2,\cdots I)$ are real polynomials
with $P_l(0)=0$, $p_1,p_2,\cdots,p_I$ runs over the prime numbers.
Let \be V(s)=Q(-\frac{1}{L}\frac{d}{ds})\zeta(s),\ee where $Q$ is a
real polynomial with $Q(0)=1$ and $Q'(x)=Q'(1-x)$. Then we have \be
\kappa\geq
1-\frac{1}{R}\log(\frac{1}{T}\int_1^T|BV(\sigma_0+it)|^2dt)+o(1).\ee

\vskip 3mm

Let $\alpha,\beta$ be complex numbers with
$\alpha,\beta\ll\frac{1}{L}$, $s_0=\frac12+iw$ with $T\leq w\leq
2T$. Let $\Delta=T^{1-\delta}$, $0<\delta<1$, $a_j$ defined by
(2.4). Let
$${\cal B}(s)=\sum_{j\leq
y}\frac{a_j}{j^s}=B(s-\frac{R}{L}),$$  and \be g(\alpha, \beta,
w)=\frac{1}{i\Delta\pi^{\frac12}}\int_{(\frac12)}e^{(s-s_0)^2\Delta^{-2}}\zeta(s+\alpha)\zeta(1-s+\beta){\cal
B}(s){\cal B}(1-s)ds,\ee where $(c)$ denotes the straight line path
from $c-i\infty$ to $c+i\infty$. By the method of Balasubramanian,
Conrey and Heath-Brown [2], to estimate the mean value integral in
(2.6), it suffice to obtain an evaluation of $g(\alpha, \beta, w)$
uniformly for $\alpha,\beta\ll\frac{1}{L}, T\leq w\leq 2T$.

The following lemma is due to Conrey [6].

\vskip 3mm

{\bf Lemma 2}. Let $\alpha,\beta$ be complex numbers with
$\alpha,\beta\ll\frac{1}{L}$, $y_1, y$, $a_j$ be as in Lemma 1,
$0<\delta<1, \Delta=T^{1-\delta}$, $g$ be as in (2.7), and \be
\Sigma(\alpha,\beta)=\sum_{h,k\leq
y}\frac{a_ha_k}{h^{1+\alpha}k^{1+\beta}}(h,k)^{1+\alpha+\beta}.\ee
Then as $T\rightarrow\infty$, \be g(\alpha, \beta,
w)=\frac{\Sigma(\beta,\alpha)-e^{-(\alpha+\beta)L}\Sigma(-\alpha,-\beta)}{\alpha+\beta}+o_{\delta}(1)\ee
uniformly in $\alpha,\beta$ and $w$.

\vskip 3mm

Lemma 2 reduce the evaluation of $g(\alpha, \beta, w)$, and
therefore the evaluation of the mean value integral in (2.6), to the
evaluation of $\Sigma(\alpha,\beta)$.

Denote \be F(j,w)=\prod_{p|j}(1-\frac{1}{p^w}),\ee we have
\begin{eqnarray}
\Sigma(\alpha,\beta)&=&\sum_{h,k\leq
y}\frac{a_ha_k}{h^{1+\alpha}k^{1+\beta}}\sum_{j|(h,k)}j^{1+\alpha+\beta}F(j,1+\alpha+\beta)\nonumber\\
&=&\sum_{j\leq y}j^{1+\alpha+\beta}F(j,1+\alpha+\beta)\sum_{j|h\atop
h\leq y}\frac{a_h}{h^{1+\alpha}}\sum_{j|k\atop k\leq
y}\frac{a_k}{k^{1+\beta}}.\end{eqnarray} Let \be
E_a(j)=\sum_{j|h\atop h\leq y}\frac{a_h}{h^{1+\alpha}}.\ee In
section 3, we evaluate $E_{\alpha}(j)$. Then we use this to evaluate
$\Sigma(\alpha,\beta)$ in section 4.

We also need the following lemmas for evaluation of $E_{\alpha}(j)$
and $\Sigma(\alpha,\beta)$.

\vskip 3mm

{\bf Lemma 3} (see Conrey [4]). Let $P$ be a real polynomial with
$P(0)=0$ and \be S=\sum_{n\leq y/j\atop
(n,j)=1}\frac{\mu(n)}{n^{1+\alpha}}P(\frac{\log y/nj}{\log y}).\ee
Then we have
\begin{eqnarray} S&=&\frac{1}{F(j,1+\alpha)}\Big(\alpha P(\frac{\log y/j}{\log y})+\frac{1}{\log y}P^{\prime}(\frac{\log y/j}{\log y})\Big)\nonumber\\
&&+O(\frac{(\log\log
y)^2F_1(j,1-2\delta)}{\log^2y})+O(\frac{(\log\log
y)^2F_1(j,1-2\delta)}{\log y}(\frac{j}{y})^d)\end{eqnarray}
uniformly for $j\leq y$, $\alpha\ll\frac{1}{\log y}$. Here $F(j,w)$
defined by (2.11), \be F_1(j,w)=\prod_{p|j}(1+\frac{1}{p^w}),\ee
$\delta=1/\log\log y$, and
 \be d=\frac{1}{M\log\log y},\ee where $M$ is a sufficiendly large constant.

\vskip 3mm

{\bf Lemma 4} (Mertens Theorem).  \be \sum_{p\leq y}\frac{\log
p}{p}=\log y+O(1).\ee

\vskip 3mm

{\bf Lemma 5} (Levinson [14]). \be \sum_{p|j}\frac{\log
p}{p}=O(\log\log j).\ee

\vskip 3mm

{\bf Lemma 6} (Levinson [16]). Let $N$ be a positive integer and \be
J(x)=\sum_{n\leq x\atop
(n,N)=1}\frac{\mu^2(n)}{n}\prod_{p|n}(1+f(p)),\ee where\be
f(p)=O(\frac{1}{p^c})\ee for some $c>0$. Then \be
J(x)=\prod_{p|N}(1-\frac{1}{p})\prod_{(p,N)=1}(1-\frac{1}{p^2})(1+\frac{f(p)}{p+1})\log
x+O(\log\log(N+1))\ee with the $O$ independent of $x$ and $N$.

\vskip 3mm

{\bf Lemma 7}. Let $m$ be a positive integer, $\alpha$ be a complex
number and $f$ be a continuous function, $D\geq 1$, then
\begin{eqnarray} &&\int_1^D\frac{1}{x_1^{1+\alpha}}\int_1^{\frac{D}{x_1}}\frac{1}{x_2^{1+\alpha}}\cdots
\int_1^{\frac{D}{x_1x_2\cdots
x_{m-2}}}\frac{1}{x_{m-1}^{1+\alpha}}\int_1^{\frac{D}{x_1x_2\cdots
x_{m-1}}}\frac{f(x_1x_2\cdots x_m)}{x_m^{1+\alpha}}
dx_m\cdots dx_2dx_1\nonumber\\
&=&\int_1^D\frac{f(x)\log^{m-1}x}{(m-1)!x^{1+\alpha}}dx.\end{eqnarray}

\vskip 3mm

{\bf Proof}.\ We prove by induction. The case $m=1$ is obvious. We
assume the case $m-1$ is valid, then by integration by parts we get
\begin{eqnarray} &&\int_1^D\frac{1}{x_1^{1+\alpha}}\int_1^{\frac{D}{x_1}}\frac{1}{x_2^{1+\alpha}}\cdots
\int_1^{\frac{D}{x_1x_2\cdots
x_{m-2}}}\frac{1}{x_{m-1}^{1+\alpha}}\int_1^{\frac{D}{x_1x_2\cdots
x_{m-1}}}\frac{f(x_1x_2\cdots x_m)}{x_m^{1+\alpha}}
dx_m\cdots dx_2dx_1\nonumber\\
&=&\int_1^D\frac{1}{x_1^{1+\alpha}}\int_1^{\frac{D}{x_1}}\frac{f(x_1x)\log^{m-2}x}{(m-2)!x^{1+\alpha}}dxdx_1
\nonumber\\
&=&\int_1^D\frac{1}{x}\int_x^D\frac{f(t)\log^{m-2}(t/x)}{(m-2)!t^{1+\alpha}}dtdx\nonumber\\
&=&\int_1^D\frac{\log
x}{x}\int_x^D\frac{f(t)\log^{m-3}(t/x)}{(m-3)!t^{1+\alpha}}dtdx\nonumber\\
&=&\cdots\nonumber\\
&=&\int_1^D\frac{\log^{m-2}x}{(m-2)!x}\int_x^D\frac{f(t)}{t^{1+\alpha}}dtdx\nonumber\\
&=&\int_1^D\frac{f(x)\log^{m-1}x}{(m-1)!x^{1+\alpha}}dx.\end{eqnarray}
The proof is complete.

\vskip 3mm

{\bf Lemma 8}. For positive integer $m_1,\ m_2$ and square-free j,
\begin{eqnarray}
&&\sum_{p_1p_2\cdots p_{m_1}|j}\log p_1\log p_2\cdots\log p_{m_1}
\sum_{q_1q_2\cdots q_{m_2}|j}\log q_1\log q_2\cdots\log
q_{m_2}\nonumber\\
&=&\sum_{k=0}^{\min(m_1,m_2)}{\cal P}_{m_1}^k{\cal
C}_{m_2}^k\sum_{p_1p_2\cdots
p_{m_1+m_2-k}|j}\log^2p_1\log^2p_2\cdots\log^2p_k\log
p_{k+1}\cdots\log p_{m_1+m_2-k},\nonumber\\&&
\end{eqnarray}
where $p$ and $q$ runs over prime numbers, ${\cal
P}_m^k=\frac{m!}{(m-k)!}$, ${\cal C}_m^k=\frac{m!}{k!(m-k)!}$.

\vskip 3mm

{\bf Proof}.\ The summation due to the case that there are just $k$
prime-square factors in $p_1p_2\cdots p_{m_1}q_1q_2\cdots q_{m_2}$
is \be {\cal P}_{m_1}^k{\cal C}_{m_2}^k\sum_{p_1p_2\cdots
p_{m_1+m_2-k}|j}\log^2p_1\log^2p_2\cdots\log^2p_k\log
p_{k+1}\cdots\log p_{m_1+m_2-k}.\ee Sum $k$ form $0$ to
$\min(m_1,m_2)$ we get (2.24).

\vskip 3mm

{\bf Lemma 9}. Let $k_1\geq 0, k_2\geq 1$ be integers, $f$ be a
continuous function, $D\geq 1$, then

\begin{eqnarray} &&\int_1^D\frac{\log x_1}{x_1}\int_1^{\frac{D}{x_1}}\frac{\log x_2}{x_2}\cdots\int_1^{\frac{D}{x_1x_2\cdots
x_{k_1-1}}}\frac{\log x_{k_1}}{x_{k_1}}\int_1^{\frac{D}{x_1x_2\cdots
x_{k_1}}}\frac{1}{x_{k_1+1}}\cdots\nonumber\\&&\int_1^{\frac{D}{x_1x_2\cdots
x_{k_1+k_2-2}}}\frac{1}{x_{k_1+k_2-1}}\int_1^{\frac{D}{x_1x_2\cdots
x_{k_1+k_2-1}}}\frac{f(x_1x_2\cdots
x_{k_1+k_2})}{x_{k_1+k_2}}dx_{k_1+k_2}\cdots dx_1
\nonumber\\
&=&\int_1^D\frac{f(x)\log^{2k_1+k_2-1}x}{(2k_1+k_2-1)!x}dx.\end{eqnarray}

\vskip 3mm

{\bf Proof}.\ We prove by induction for $k_1$. Let $\alpha=0$ and
$m=k_2$ in lemma 7, the case $k_1=0$ follows. We assume the case
$k_1-1$ is valid, then similar to (2.23), by integration by parts we
get
\begin{eqnarray}&&\int_1^D\frac{\log x_1}{x_1}\int_1^{\frac{D}{x_1}}\frac{\log x_2}{x_2}\cdots\int_1^{\frac{D}{x_1x_2\cdots
x_{k_1-1}}}\frac{\log x_{k_1}}{x_{k_1}}\int_1^{\frac{D}{x_1x_2\cdots
x_{k_1}}}\frac{1}{x_{k_1+1}}\cdots\nonumber\\&&\int_1^{\frac{D}{x_1x_2\cdots
x_{k_1+k_2-2}}}\frac{1}{x_{k_1+k_2-1}}\int_1^{\frac{D}{x_1x_2\cdots
x_{k_1+k_2-1}}}\frac{f(x_1x_2\cdots
x_{k_1+k_2})}{x_{k_1+k_2}}dx_{k_1+k_2}\cdots dx_1
\nonumber\\
&=&\int_1^D\frac{\log
x_1}{x_1}\int_1^{\frac{D}{x_1}}\frac{f(x_1x)\log^{2k_1+k_2-3}x}{(2k_1+k_2-3)!x}dxdx_1
\nonumber\\
&=&\int_1^D\frac{\log
x}{x}\int_x^D\frac{f(t)\log^{2k_1+k_2-3}(t/x)}{(2k_1+k_2-3)!t}dtdx\nonumber\\
&=&\int_1^D\frac{\log^2
x}{2x}\int_x^D\frac{f(t)\log^{2k_1+k_2-4}(t/x)}{(2k_1+k_2-4)!t}dtdx\nonumber\\
&=&\cdots\nonumber\\
&=&\int_1^D\frac{\log^{2k_1+k_2-2}x}{(2k_1+k_2-2)!x}\int_x^D\frac{f(t)}{t}dtdx
\nonumber\\
&=&\int_1^D\frac{f(x)\log^{2k_1+k_2-1}x}{(2k_1+k_2-1)!x}dx.\end{eqnarray}

 \vskip 8mm
\def\theequation{3.\arabic{equation}}
\setcounter{equation}{0} \centerline{\large\bf 3.\ Evaluation of
$E_{\alpha}(j)$ } \vskip 5mm

\noindent Throughout this section, estimation are uniformly for
$j\leq y$, $\alpha\ll\frac{1}{\log y}$.

By (2.12) and (2.4), we have
\begin{eqnarray}  E_{\alpha}(j)&=&\sum_{j|h\atop
h\leq y}\frac{\mu(h)}{h^{1+\alpha}}P_1(\frac{\log y/h}{\log
y})\nonumber\\
&&+\sum_{j|h\atop h\leq
y_1}\frac{\mu(h)}{h^{1+\alpha}}\sum_{p_1p_2|h}\frac{\log p_1\log
p_2}{\log^2y_1}P_2(\frac{\log
y_1/h}{\log y_1})+\cdots\nonumber\\
&&+\sum_{j|h\atop h\leq
y_1}\frac{\mu(h)}{h^{1+\alpha}}\sum_{p_1p_2\cdots p_I|h}\frac{\log
p_1\log p_2\cdots\log p_I}{\log^Iy_1}P_I(\frac{\log
y_1/h}{\log y_1})\nonumber\\
&=&\Sigma_1+\Sigma_2+\cdots+\Sigma_I,\end{eqnarray} say.

Let $n=h/j$, $\delta=1/\log\log y$, and $ d=\frac{1}{M\log\log y}$,
where $M$ is a sufficiendly large constant, by lemma 3 we have
\begin{eqnarray} \Sigma_1&=&\frac{\mu(j)}{j^{1+\alpha}}\sum_{n\leq y/j\atop (n,j)=1}\frac{\mu(n)}{n^{1+\alpha}}P_1(\frac{\log (y/nj)}{\log
y})\nonumber\\
&=&\frac{\mu(j)}{j^{1+\alpha}F(j,1+\alpha)}\Big(\alpha P_1(\frac{\log y/j}{\log y})+\frac{1}{\log y}P_1^{\prime}(\frac{\log y/j}{\log y})\Big)\nonumber\\
&&+O\Big(\frac{\mu(j)(\log\log
y)^2F_1(j,1-2\delta)}{j\log^2y}\Big)\nonumber\\
&&+O\Big(\frac{\mu(j)(\log\log y)^2F_1(j,1-2\delta)}{j\log
y}(\frac{j}{y})^d\Big).\end{eqnarray} For $l\geq 2$, we only need
consider $j\leq y_1$, since $\Sigma_l=0$ for $j>y_1$. Let $n=h/j$,
we have
\begin{eqnarray} \Sigma_l&=&\frac{\mu(j)}{j^{1+\alpha}}\sum_{n\leq y_1/j\atop (n,j)=1}\frac{\mu(n)}{n^{1+\alpha}}\sum_{p_1p_2\cdots
p_l|nj}\frac{\log p_1\log p_2\cdots\log
p_l}{\log^ly_1}P_l(\frac{\log (y_1/nj)}{\log
y_1})\nonumber\\
&=&\frac{\mu(j)}{j^{1+\alpha}}\sum_{n\leq y_1/j\atop
(n,j)=1}\frac{\mu(n)}{n^{1+\alpha}}P_l(\frac{\log (y_1/nj)}{\log
y_1})\Big(\sum_{p_1p_2\cdots p_l|j}\frac{\log p_1\log p_2\cdots\log
p_l}{\log^ly_1}\nonumber\\
&&\hspace*{8mm}+\cdots
\nonumber\\
&&\hspace*{8mm}+{\cal C}_l^m\sum_{p_1p_2\cdots p_{l-m}|j}\frac{\log
p_1\log p_2\cdots\log p_{l-m}}{\log^{l-m}y_1}\sum_{p_1p_2\cdots
p_m|n}\frac{\log p_1\log
p_2\cdots \log p_m}{\log^m y_1}\nonumber\\
&&\hspace*{8mm}+\cdots\nonumber\\
&&\hspace*{8mm}+\sum_{p_1p_2\cdots p_l|n}\frac{\log p_1\log
p_2\cdots\log p_l}{\log^ly_1}\Big).\end{eqnarray}
 For $1\leq m\leq l$, let
$n_0=\frac{n}{p_1p_2\cdots p_m}$, we have by lemma 3
\begin{eqnarray}&&\sum_{n\leq y_1/j\atop
(n,j)=1}\frac{\mu(n)}{n^{1+\alpha}}P_l(\frac{\log (y_1/nj)}{\log
y_1})\sum_{p_1p_2\cdots p_m|n}\frac{\log p_1\log p_2\cdots\log
p_m}{\log^my_1}\nonumber\\
&=&\sum_{p_1p_2\cdots p_m\leq y_1/j\atop(p_1p_2\cdots
p_m,j)=1}\frac{\mu(p_1p_2\cdots p_m)\log p_1\log p_2\cdots\log
p_m}{(p_1p_2\cdots
p_m)^{1+\alpha}\log^my_1}\nonumber\\
&&\hspace*{4cm}\times\sum_{n_0\leq \frac{y_1}{p_1p_2\cdots
p_mj}\atop(n_0,p_1p_2\cdots
p_mj)=1}\frac{\mu(n_0)}{n_0^{1+\alpha}}P_l(\frac{\log
(y_1/n_0p_1p_2\cdots p_mj)}{\log y_1})\nonumber\\
&=&\sum_{p_1p_2\cdots p_m\leq y_1/j\atop(p_1p_2\cdots
p_m,j)=1}\frac{\mu(p_1p_2\cdots p_m)\log p_1\log p_2\cdots\log
p_m}{F(p_1p_2\cdots p_mj,1+\alpha)(p_1p_2\cdots
p_m)^{1+\alpha}\log^my_1}\nonumber\\
&&\hspace*{4cm}\times(\alpha P_l(\frac{\log\frac{y_1}{p_1p_2\cdots
p_mj}}{\log y_1})+\frac{1}{\log
y_1}P_l^{\prime}(\frac{\log\frac{y_1}{p_1p_2\cdots p_mj}}{\log
y_1}))
\nonumber\\
&&+O\Big(\sum_{p_1p_2\cdots p_m\leq y_1/j\atop(p_1p_2\cdots
p_m,j)=1}\frac{\log p_1\log p_2\cdots\log p_m(\log\log
y_1)^2F_1(p_1p_2\cdots p_mj,1-2\delta_1)}{p_1p_2\cdots
p_m\log^{m+2}y_1}\Big)
\nonumber\\
&&+O\Big(\sum_{p_1p_2\cdots p_m\leq y_1/j\atop(p_1p_2\cdots
p_m,j)=1}\frac{\log p_1\log p_2\cdots\log p_m(\log\log
y_1)^2F_1(p_1p_2\cdots p_mj,1-2\delta_1)}{p_1p_2\cdots
p_m\log^{m+1}y_1}\nonumber\\
&&\hspace*{3cm}\times(\frac{p_1p_2\cdots
p_mj}{y_1})^{d_1}\Big)\nonumber\\
&=&A_1+O(A_2)+O(A_3),\end{eqnarray} say, where $\delta_1=1/\log\log
y_1$, $d_1=\frac{1}{M\log\log y_1}$ with $M$ a sufficiendly large
constant.

Now we evaluate $A_1$ first. We have by (2.10)
\begin{eqnarray}&&\sum_{p_1p_2\cdots p_m\leq y_1/j\atop(p_1p_2\cdots
p_m,j)=1}\frac{\mu(p_1p_2\cdots p_m)\log p_1\log p_2\cdots\log
p_m}{F(p_1p_2\cdots p_mj,1+\alpha)(p_1p_2\cdots
p_m)^{1+\alpha}}P_l(\frac{\log\frac{y_1}{p_1p_2\cdots
p_mj}}{\log y_1})\nonumber\\
&=&\frac{-1}{F(j,1+\alpha)}\sum_{p_1p_2\cdots p_{m-1}\leq
y_1/j\atop(p_1p_2\cdots p_{m-1},j)=1}\frac{\mu(p_1p_2\cdots
p_{m-1})\log p_1\log p_2\cdots\log p_{m-1}}{(p_1^{1+\alpha}-1)(p_2^{1+\alpha}-1)\cdots(p_{m-1}^{1+\alpha}-1)}\nonumber\\
&&\hspace*{3cm}\times\sum_{p_m\leq y_1/p_1p_2\cdots p_{m-1}j\atop
(p_m, p_1p_2\cdots p_{m-1}j)=1}\frac{\log
p_m}{p_m^{1+\alpha}-1}P_l(\frac{\log\frac{y_1}{p_1p_2\cdots
p_mj}}{\log y_1}).
\end{eqnarray}
Consider the inner summation, we have
\begin{eqnarray}&&\sum_{p_m\leq y_1/p_1p_2\cdots p_{m-1}j\atop (p_m,
p_1p_2\cdots p_{m-1}j)=1}\frac{\log
p_m}{p_m^{1+\alpha}-1}P_l(\frac{\log\frac{y_1}{p_1p_2\cdots
p_mj}}{\log
y_1})\nonumber\\
&=&\sum_{p_m\leq y_1/p_1p_2\cdots p_{m-1}j}\frac{\log
p_m}{p_m^{1+\alpha}-1}P_l(\frac{\log\frac{y_1}{p_1p_2\cdots
p_mj}}{\log
y_1})\nonumber\\
&&-\sum_{p_m\leq y_1/p_1p_2\cdots p_{m-1}j\atop p_m|p_1p_2\cdots
p_{m-1}j}\frac{\log
p_m}{p_m^{1+\alpha}-1}P_l(\frac{\log\frac{y_1}{p_1p_2\cdots
p_mj}}{\log
y_1})\nonumber\\
&=&B_1-B_2,\end{eqnarray} say. Since $P_l(x)$ is bounded in $[0,1]$,
we have for $p_1p_2\cdots p_mj\leq y_1$, \be
P_l(\frac{\log\frac{y_1}{p_1p_2\cdots p_mj}}{\log y_1})=O(1).\ee
Thus by lemma 4 and Abel summation,
\begin{eqnarray}B_1&=&\int_1^{\frac{y_1}{p_1p_2\cdots p_{m-1}j}}\frac{1}{x_m^{1+\alpha}-1}P_l(\frac{\log\frac{y_1}{p_1p_2\cdots
p_{m-1}jx_m}}{\log
y_1})dx_m+O(1)\nonumber\\
&=&\int_1^{\frac{y_1}{p_1p_2\cdots
p_{m-1}j}}\frac{1}{x_m^{1+\alpha}}P_l(\frac{\log\frac{y_1}{p_1p_2\cdots
p_{m-1}jx_m}}{\log y_1})dx_m+O(1),\end{eqnarray} by lemma 5, \be
B_2=O(\log\log( p_1p_2\cdots p_{m-1}j))=O(\log\log y_1).\ee Combine
(3.5), (3.6), (3.8) and (3.9), we get
\begin{eqnarray}&&\sum_{p_1p_2\cdots p_m\leq y_1/j\atop(p_1p_2\cdots
p_m,j)=1}\frac{\mu(p_1p_2\cdots p_m)\log p_1\log p_2\cdots\log
p_m}{F(p_1p_2\cdots p_mj,1+\alpha)(p_1p_2\cdots
p_m)^{1+\alpha}}P_l(\frac{\log\frac{y_1}{p_1p_2\cdots
p_mj}}{\log y_1})\nonumber\\
&=&\frac{-1}{F(j,1+\alpha)}\sum_{p_1p_2\cdots p_{m-1}\leq
y_1/j\atop(p_1p_2\cdots p_{m-1},j)=1}\frac{\mu(p_1p_2\cdots
p_{m-1})\log p_1\log p_2\cdots\log p_{m-1}}{(p_1^{1+\alpha}-1)(p_2^{1+\alpha}-1)\cdots(p_{m-1}^{1+\alpha}-1)}\nonumber\\
&&\hspace*{5cm}\times \int_1^{\frac{y_1}{p_1p_2\cdots
p_{m-1}j}}\frac{1}{x_m^{1+\alpha}}P_l(\frac{\log\frac{y_1}{p_1p_2\cdots
p_{m-1}jx_m}}{\log y_1})dx_m\nonumber\\
&&+O\Huge(\frac{\log\log y_1}{F(j,1+\alpha)}\sum_{p_1p_2\cdots
p_{m-1}\leq y_1/j\atop(p_1p_2\cdots
p_{m-1},j)=1}\frac{\mu(p_1p_2\cdots p_{m-1})\log p_1\log
p_2\cdots\log
p_{m-1}}{(p_1^{1+\alpha}-1)(p_2^{1+\alpha}-1)\cdots(p_{m-1}^{1+\alpha}-1)}\Huge)\nonumber\\&=&C_1+O(C_2),
\end{eqnarray} say. By lemma 4, we have
\begin{eqnarray} C_2&=&O\Huge(F_1(j,1-2\delta_1)\log\log y_1\prod_{r=1}^{m-1}\sum_{p_r\leq y_1}|\frac{\log p_r}{p_r^{1+\alpha}-1}|\huge)
\nonumber\\
&=&O(F_1(j,1-2\delta_1)\log^{m-1}y_1\log\log y_1).\end{eqnarray}
(3.10) and (3.11) gives
\begin{eqnarray}&&\sum_{p_1p_2\cdots p_m\leq y_1/j\atop(p_1p_2\cdots
p_m,j)=1}\frac{\mu(p_1p_2\cdots p_m)\log p_1\log p_2\cdots\log
p_m}{F(p_1p_2\cdots p_mj,1+\alpha)(p_1p_2\cdots
p_m)^{1+\alpha}}P_l(\frac{\log\frac{y_1}{p_1p_2\cdots
p_mj}}{\log y_1})\nonumber\\
&=&\frac{1}{F(j,1+\alpha)}\sum_{p_1p_2\cdots p_{m-2}\leq
y_1/j\atop(p_1p_2\cdots p_{m-2},j)=1}\frac{\mu(p_1p_2\cdots
p_{m-2})\log p_1\log p_2\cdots\log p_{m-2}}{(p_1^{1+\alpha}-1)(p_2^{1+\alpha}-1)\cdots(p_{m-2}^{1+\alpha}-1)}\nonumber\\
&&\hspace*{5mm}\times\sum_{p_{m-1}\leq y_1/p_1p_2\cdots
p_{m-2}j\atop (p_{m-1}, p_1p_2\cdots p_{m-2}j)=1}\frac{\log
p_{m-1}}{p_{m-1}^{1+\alpha}-1}\int_1^{\frac{y_1}{p_1p_2\cdots
p_{m-1}j}}\frac{1}{x_m^{1+\alpha}}P_l(\frac{\log\frac{y_1}{p_1p_2\cdots
p_{m-1}jx_m}}{\log y_1})dx_m\nonumber\\
&&+O(F_1(j,1-2\delta_1)\log^{m-1}y_1\log\log y_1).
\end{eqnarray}
Therefore by induction for $m$ and lemma 7, we obtain
\begin{eqnarray}
&&\sum_{p_1p_2\cdots p_m\leq y_1/j\atop(p_1p_2\cdots
p_m,j)=1}\frac{\mu(p_1p_2\cdots p_m)\log p_1\log p_2\cdots\log
p_m}{F(p_1p_2\cdots p_mj,1+\alpha)(p_1p_2\cdots
p_m)^{1+\alpha}}P_l(\frac{\log\frac{y_1}{p_1p_2\cdots
p_mj}}{\log y_1})\nonumber\\
&=&\frac{(-1)^m}{F(j,1+\alpha)}\int_1^{\frac{y_1}{j}}\frac{1}{x_1^{1+\alpha}}\cdots
\int_1^{\frac{y_1}{jx_1\cdots
x_{m-2}}}\frac{1}{x_{m-1}^{1+\alpha}}\nonumber\\
&&\hspace*{1cm}\times\int_1^{\frac{y_1}{jx_1\cdots
x_{m-1}}}\frac{1}{x_m^{1+\alpha}}
P_l(\frac{\log\frac{y_1}{jx_1\cdots
x_{m-1}x_m}}{\log y_1})dx_mdx_{m-1}\cdots dx_1\nonumber\\
&&+O(F_1(j,1-2\delta_1)\log^{m-1}y_1\log\log y_1)\nonumber\\
&=&\frac{(-1)^m}{F(j,1+\alpha)(m-1)!}\int_1^{\frac{y_1}{j}}\frac{P_l(\frac{\log
y_1/jx}{\log
y_1})\log^{m-1}x}{x^{1+\alpha}}dx\nonumber\\
&&+O(F_1(j,1-2\delta_1)\log^{m-1}y_1\log\log y_1).
\end{eqnarray}
Similarly,
\begin{eqnarray}
&&\sum_{p_1p_2\cdots p_m\leq y_1/j\atop(p_1p_2\cdots
p_m,j)=1}\frac{\mu(p_1p_2\cdots p_m)\log p_1\log p_2\cdots\log
p_m}{F(p_1p_2\cdots p_mj,1+\alpha)(p_1p_2\cdots
p_m)^{1+\alpha}}P_l^{\prime}(\frac{\log\frac{y_1}{p_1p_2\cdots
p_mj}}{\log y_1})\nonumber\\
&=&\frac{(-1)^m}{F(j,1+\alpha)(m-1)!}\int_1^{\frac{y_1}{j}}\frac{P_l^{\prime}(\frac{\log
y_1/jx}{\log
y_1})\log^{m-1}x}{x^{1+\alpha}}dx\nonumber\\
&&+O(F_1(j,1-2\delta_1)\log^{m-1}y_1\log\log y_1).
\end{eqnarray}
Since $P_l(0)=0$, we have \begin{eqnarray}
&&\int_1^{\frac{y_1}{j}}\frac{(\alpha P_l(\frac{\log y_1/jx}{\log
y_1})+\frac{1}{\log y_1}P_l^{\prime}(\frac{\log y_1/jx}{\log
y_1}))}{x^{1+\alpha}}dx\nonumber\\
&=&\int_1^{\frac{y_1}{j}}\alpha P_l(\frac{\log y_1/jx}{\log y_1})
\frac{1}{x^{1+\alpha}}dx-P_l(\frac{\log y_1/jx}{\log
y_1})\frac{1}{x^{\alpha}}\Big|_1^{y_1/j}+\int_1^{\frac{y_1}{j}}P_l(\frac{\log
y_1/jx}{\log y_1})(\frac{1}{x^{\alpha}})^{\prime}dx\nonumber\\
&=&P_l(\frac{\log y_1/j}{\log y_1}),\end{eqnarray} and for $m\geq
2$,
\begin{eqnarray}
&&\int_1^{\frac{y_1}{j}}\frac{(\alpha P_l(\frac{\log y_1/jx}{\log
y_1})+\frac{1}{\log y_1}P_l^{\prime}(\frac{\log y_1/jx}{\log
y_1}))\log^{m-1}x}{x^{1+\alpha}}dx\nonumber\\
&=&\int_1^{\frac{y_1}{j}}\alpha P_l(\frac{\log y_1/jx}{\log y_1})
\frac{\log^{m-1}x}{x^{1+\alpha}}dx\nonumber\\
&&-P_l(\frac{\log y_1/jx}{\log
y_1})\frac{\log^{m-1}x}{x^{\alpha}}\Big|_1^{y_1/j}+\int_1^{\frac{y_1}{j}}P_l(\frac{\log
y_1/jx}{\log y_1})(\frac{\log^{m-1}x}{x^{\alpha}})^{\prime}dx\nonumber\\
&=&(m-1)\int_1^{\frac{y_1}{j}}P_l(\frac{\log y_1/jx}{\log
y_1})\frac{\log^{m-2}x}{x^{1+\alpha}}dx.
\end{eqnarray}
Combine (3.13)-(3.16), we have for $m=1$, \be
A_1=\frac{-1}{F(j,1+\alpha)\log y_1}P_l(\frac{\log y_1/j}{\log
y_1})+O(\frac{F_1(j,1-2\delta_1)\log\log y_1}{\log^2 y_1})\ee and
for $m\geq 2$,
\begin{eqnarray}
A_1&=&\frac{(-1)^m}{F(j,1+\alpha)(m-2)!\log^my_1}\int_1^{\frac{y_1}{j}}P_l(\frac{\log
y_1/jx}{\log
y_1})\frac{\log^{m-2}x}{x^{1+\alpha}}dx\nonumber\\
&&+O(\frac{F_1(j,1-2\delta_1)\log\log y_1}{\log^2 y_1}).
\end{eqnarray}

For $(p_1p_2\cdots p_m,j)=1$, we have \begin{eqnarray}&&
F_1(p_1p_2\cdots p_mj,1-2\delta_1)=F_1(j,
1-2\delta_1)F_1(p_1p_2\cdots
p_m,1-2\delta_1)\nonumber\\
&\leq &2^mF_1(j, 1-2\delta_1),\end{eqnarray} hence by lemma 4 we
obtain
\begin{eqnarray} A_2&=&O\Huge(\frac{F_1(j,1-2\delta_1)(\log\log y_1)^2}{\log^{m+2}y_1}\prod_{r=1}^m\sum_{p_r\leq y_1}
\frac{\log p_r}{p_r}\Huge)
\nonumber\\
&=&O\Huge(\frac{F_1(j,1-2\delta_1)(\log\log
y_1)^2}{\log^2y_1}\Huge).\end{eqnarray} By lemma 4 and Abel
summation we have \be \sum_{p_m\leq y_1/p_1p_2\cdots
p_{m-1}j}\frac{\log p_m}{p_m^{1-d_1}}=O(\frac{y_1^{d_1}\log\log
y_1}{(p_1p_2\cdots p_{m-1}j)^{d_1}}).\ee (3.19) and (3.21) yields
\begin{eqnarray}A_3&=&O\Big(\frac{F(j,1-2\delta_1)(\log\log y_1)^2}{\log^{m+1}y_1}\sum_{p_1p_2\cdots p_{m-1}\leq y_1/j}\frac{\log p_1\log p_2\cdots\log p_{m-1}}{p_1p_2\cdots
p_{m-1}}\nonumber\\
&&\hspace*{4cm}\times\frac{(p_1p_2\cdots
p_{m-1}j)^{d_1}}{y^{d_1}}\sum_{p_m\leq y_1/p_1p_2\cdots
p_{m-1}j}\frac{\log
p_m}{p_m^{1-d_1}}\Big)\nonumber\\
&=&O\Big(\frac{F(j,1-2\delta_1)(\log\log
y_1)^3}{\log^{m+1}y_1}\prod_{r=1}^{m-1}\sum_{p_r\leq y_1}\frac{\log p_r}{p_r}\Big)\nonumber\\
&=&O\Big(\frac{F(j,1-2\delta_1)(\log\log y_1)^3}{\log^2y_1}\Big).
\end{eqnarray}

Substitute (3.17), (3.18), (3.20) and (3.22) into (3.4), we obtain
\begin{eqnarray} &&\sum_{n\leq y_1/j\atop
(n,j)=1}\frac{\mu(n)}{n^{1+\alpha}}P_l(\frac{\log (y_1/nj)}{\log
y_1})\sum_{p_1|n}\frac{\log p_1}{\log
y_1}\nonumber\\
&=&\frac{-1}{F(j,1+\alpha)\log y_1}P_l(\frac{\log y_1/j}{\log
y_1})+O(\frac{F_1(j,1-2\delta_1)(\log\log y_1)^3}{\log^2
y_1})\end{eqnarray} and for $m\geq 2$,
\begin{eqnarray} &&\sum_{n\leq y_1/j\atop
(n,j)=1}\frac{\mu(n)}{n^{1+\alpha}}P_l(\frac{\log (y_1/nj)}{\log
y_1})\sum_{p_1p_2\cdots p_m|n}\frac{\log p_1\log p_2\cdots\log
p_m}{\log^my_1}\nonumber\\
&=&\frac{(-1)^m}{F(j,1+\alpha)(m-2)!\log^my_1}\int_1^{\frac{y_1}{j}}P_l(\frac{\log
y_1/jx}{\log
y_1})\frac{\log^{m-2}x}{x^{1+\alpha}}dx\nonumber\\
&&+O(\frac{F_1(j,1-2\delta_1)(\log\log y_1)^3}{\log^2
y_1}).\end{eqnarray} For $j\leq y_1$, we have trivially \be
\sum_{p_1p_2\cdots p_m|j}\frac{\log p_1\log p_2\cdots\log
p_m}{\log^m y_1}\leq \prod_{r=1}^m\sum_{p_r|j}\frac{\log p_r}{\log
y_1}\leq 1.\ee Substitute (3.23)-(3.25) and (3.2) with $P_1, y,
\delta,d$ replaced by $P_l, y_1,\delta_1,d_1$ into (3.3), we get for
$j\leq y_1$,
\begin{eqnarray}&&\Sigma_l=\frac{\mu(j)}{j^{1+\alpha}F(j,1+\alpha)\log^ly_1}\Big\{\frac{(-1)^l}{(l-2)!}
\int_1^{\frac{y_1}{j}}P_l(\frac{\log y_1/jx}{\log
y_1})\frac{\log^{l-2}x}{x^{1+\alpha}}dx\nonumber\\
&&+\sum_{m=1}^{l-2}{\cal
C}_l^m\frac{(-1)^{l-m}}{(l-m-2)!}\int_1^{\frac{y_1}{j}}P_l(\frac{\log
y_1/jx}{\log
y_1})\frac{\log^{l-m-2}x}{x^{1+\alpha}}dx\sum_{p_1p_2\cdots
p_m|j}\log p_1\log
p_2\cdots \log p_m\nonumber\\
&&-{\cal C}_l^{l-1}P_l(\frac{\log y_1/j}{\log
y_1})\sum_{p_1p_2\cdots p_{l-1}|j}\log p_1\log p_2\cdots\log
p_{l-1}\nonumber\\
&&+\Big(\alpha P_l(\frac{\log y_1/j}{\log y_1})+\frac{1}{\log
y_1}P_l^{\prime}(\frac{\log y_1/j}{\log y_1})\Big)\sum_{p_1p_2\cdots
p_l|j}\log p_1\log p_2\cdots\log
p_l\Big\}\nonumber\\
&&+O(\frac{\mu(j)F_1(j,1-2\delta_1)(\log\log y_1)^3}{j\log^2
y_1})+O(\frac{\mu(j)F_1(j,1-2\delta_1)(\log\log y_1)^2}{j\log
y_1}(\frac{j}{y_1})^{d_1}).\end{eqnarray}

Let $1_{y_1}(j)=1$ for $ j\leq y_1$, and $1_{y_1}(j)=0$ for $
j>y_1$, substitute (3.26) and (3.2) into (3.1), we have for $j\leq
y$,
\begin{eqnarray}
&&E_{\alpha}(j)=\frac{\mu(j)}{j^{1+\alpha}F(j,1+\alpha)}\Big\{G_0(\alpha,j)+G_1(\alpha,j)\sum_{p_1|j}\log
p_1+G_2(\alpha,j)\sum_{p_1p_2|j}\log p_1\log p_2\nonumber\\
&&\hspace*{4cm}+\cdots+G_I(\alpha,j)\sum_{p_1p_2\cdots p_I|j}\log
p_1\log
p_2\cdots\log p_I\Big\}\nonumber\\
&&\hspace*{5mm}+O(\frac{\mu(j)F_1(j,1-2\delta)(\log\log
y)^3}{j\log^2 y})+O(\frac{\mu(j)F_1(j,1-2\delta)(\log\log
y)^2}{j\log y}(\frac{j}{y})^d)\nonumber\\
&&\hspace*{5mm}+O(1_{y_1}(j)\frac{\mu(j)F_1(j,1-2\delta_1)(\log\log
y_1)^3}{j\log^2 y_1})\nonumber\\
&&\hspace*{5mm}+O(1_{y_1}(j)\frac{\mu(j)F_1(j,1-2\delta_1)(\log\log
y_1)^2}{j\log y_1}(\frac{j}{y_1})^{d_1}),\end{eqnarray} where
\begin{eqnarray}  G_0(\alpha,j)&=&
\alpha P_1(\frac{\log y/j}{\log y})+\frac{1}{\log
y}P_1^{\prime}(\frac{\log y/j}{\log
y})\nonumber\\&&+1_{y_1}(j)\sum_{l=2}^I\frac{(-1)^l}{(l-2)!\log^ly_1}
\int_1^{\frac{y_1}{j}}P_l(\frac{\log y_1/jx}{\log
y_1})\frac{\log^{l-2}x}{x^{1+\alpha}}dx,\end{eqnarray}
\begin{eqnarray}
G_1(\alpha,j)&=&1_{y_1}(j)\Big\{\frac{-2}{\log^2y_1}P_2(\frac{\log
y_1/j}{\log y_1})\nonumber\\
&&+\sum_{l=3}^I{\cal C}_l^1\frac{(-1)^{l-1}}{(l-3)!\log^ly_1}
\int_1^{\frac{y_1}{j}}P_l(\frac{\log y_1/jx}{\log
y_1})\frac{\log^{l-3}x}{x^{1+\alpha}}dx\Big\},\end{eqnarray} for
$2\leq m\leq I-2$,
\begin{eqnarray}
G_m(\alpha,j)&&=1_{y_1}(j)\Big\{\frac{1}{\log^my_1}\Big(\alpha
P_m(\frac{\log y_1/j}{\log y_1})+\frac{1}{\log
y_1}P_m^{\prime}(\frac{\log y_1/j}{\log y_1})\Big)\nonumber\\
&&+{\cal C}_{m+1}^m\frac{-1}{\log^{m+1}y_1}P_{m+1}(\frac{\log
y_1/j}{\log
y_1})\nonumber\\
&&+\sum_{l=m+2}^I{\cal C}_l^m\frac{(-1)^{l-m}}{(l-m-2)!\log^ly_1}
\int_1^{\frac{y_1}{j}}P_l(\frac{\log y_1/jx}{\log
y_1})\frac{\log^{l-m-2}x}{x^{1+\alpha}}dx\Big\},\end{eqnarray} and
\be G_{I-1}(\alpha,j)=1_{y_1}(j)\Big\{\frac{\alpha
P_{I-1}(\frac{\log y_1/j}{\log y_1})+\frac{1}{\log
y_1}P_{I-1}^{\prime}(\frac{\log y_1/j}{\log
y_1})}{\log^{I-1}y_1}-\frac{I}{\log^Iy_1}P_I(\frac{\log y_1/j}{\log
y_1})\Big\},\ee \be
G_I(\alpha,j)=1_{y_1}(j)\frac{1}{\log^Iy_1}\Big\{\alpha
P_I(\frac{\log y_1/j}{\log y_1})+\frac{1}{\log
y_1}P_I^{\prime}(\frac{\log y_1/j}{\log y_1})\Big\}.\ee

\vskip 8mm
\def\theequation{4.\arabic{equation}}
\setcounter{equation}{0} \centerline{\large\bf 4.\ Evaluation of
$\Sigma(\alpha,\beta)$ } \vskip 5mm

\noindent Throughout this section, estimation are uniformly for
$\alpha,\ \beta\ll\frac{1}{\log y}$

By (2.11) and (3.27), we have \be
\Sigma(\alpha,\beta)=U_1+U_2+\cdots+U_8+U_2'+\cdots+U_8',\ee where
\begin{eqnarray}
&&U_1=\sum_{j\leq
y}\frac{\mu^2(j)F(j,1+\alpha+\beta)}{jF(j,1+\alpha)F(j,1+\beta)}\Big(G_0(\alpha,j)+G_1(\alpha,j)\sum_{p_1|j}\log
p_1\nonumber\\
&&\hspace*{4cm}+\cdots+G_I(\alpha,j)\sum_{p_1p_2\cdots p_I|j}\log
p_1\log
p_2\cdots\log p_I\Big)\nonumber\\
&&\hspace*{2cm}\times\Big(G_0(\beta,j)+G_1(\beta,j)\sum_{p_1|j}\log
p_1\nonumber\\&&\hspace*{4cm}+\cdots+G_I(\beta,j)\sum_{p_1p_2\cdots
p_I|j}\log p_1\log p_2\cdots\log p_I\Big),\end{eqnarray} \be
U_2=O\Big(\frac{(\log\log y)^6}{\log^4y}\sum_{j\leq
y}\Big|\frac{\mu^2(j)F(j,1+\alpha+\beta)F_1^2(j,1-2\delta)}{j^{1-\alpha-\beta}}\Big|\Big),\ee
\be U_3=O\Big(\frac{(\log\log y)^4}{y^{2d}\log^2y}\sum_{j\leq
y}\Big|\frac{\mu^2(j)F(j,1+\alpha+\beta)F_1^2(j,1-2\delta)}{j^{1-2d-\alpha-\beta}}\Big|\Big),\ee
\begin{eqnarray}
U_4&=&O\Big(\frac{(\log\log y)^3}{\log^2y}\sum_{j\leq
y}\Big|\frac{\mu^2(j)F(j,1+\alpha+\beta)F_1(j,1-2\delta)}{j^{1-\beta}F(j,1+\alpha)}\nonumber\\
&&\hspace*{1cm}\times\Big(G_0(\alpha,j)+G_1(\alpha,j)\sum_{p_1|j}\log
p_1+G_2(\alpha,j)\sum_{p_1p_2|j}\log p_1\log
p_2\nonumber\\
&&\hspace*{2cm}+\cdots+G_I(\alpha,j)\sum_{p_1p_2\cdots p_I|j}\log
p_1\log p_2\cdots\log p_I\Big)\Big|\Big),\end{eqnarray}
\begin{eqnarray}
U_5&=&O\Big(\frac{(\log\log y)^2}{y^d\log y}\sum_{j\leq
y}\Big|\frac{\mu^2(j)F(j,1+\alpha+\beta)F_1(j,1-2\delta)}{j^{1-d-\beta}F(j,1+\alpha)}\nonumber\\
&&\hspace*{1cm}\times\Big(G_0(\alpha,j)+G_1(\alpha,j)\sum_{p_1|j}\log
p_1+G_2(\alpha,j)\sum_{p_1p_2|j}\log p_1\log
p_2\nonumber\\
&&\hspace*{2cm}+\cdots+G_I(\alpha,j)\sum_{p_1p_2\cdots p_I|j}\log
p_1\log p_2\cdots\log p_I\Big)\Big|\Big),\end{eqnarray}
\begin{eqnarray}
U_6&=&O\Big(\frac{(\log\log y)^3}{\log^2y}\sum_{j\leq
y}\Big|\frac{\mu^2(j)F(j,1+\alpha+\beta)F_1(j,1-2\delta)}{j^{1-\alpha}F(j,1+\beta)}\nonumber\\
&&\hspace*{1cm}\times\Big(G_0(\beta,j)+G_1(\beta,j)\sum_{p_1|j}\log
p_1+G_2(\beta,j)\sum_{p_1p_2|j}\log p_1\log
p_2\nonumber\\
&&\hspace*{2cm}+\cdots+G_I(\beta,j)\sum_{p_1p_2\cdots p_I|j}\log
p_1\log p_2\cdots\log p_I\Big)\Big|\Big),\end{eqnarray}
\begin{eqnarray}
U_7&=&O\Big(\frac{(\log\log y)^2}{y^d\log y}\sum_{j\leq
y}\Big|\frac{\mu^2(j)F(j,1+\alpha+\beta)F_1(j,1-2\delta)}{j^{1-d-\alpha}F(j,1+\beta)}\nonumber\\
&&\hspace*{1cm}\times\Big(G_0(\beta,j)+G_1(\beta,j)\sum_{p_1|j}\log
p_1+G_2(\beta,j)\sum_{p_1p_2|j}\log p_1\log
p_2\nonumber\\
&&\hspace*{2cm}+\cdots+G_I(\beta,j)\sum_{p_1p_2\cdots p_I|j}\log
p_1\log p_2\cdots\log p_I\Big)\Big|\Big),\end{eqnarray} \be U_8=
O\Big(\frac{(\log\log y)^5}{y^d\log^3y}\sum_{j\leq
y}\Big|\frac{\mu^2(j)F(j,1+\alpha+\beta)F_1^2(j,1-2\delta)}{j^{1-d-\alpha-\beta}}\Big|\Big),\ee
$U_2',U_3',\cdots, U_8'$ are the same as $U_2,U_3,\cdots, U_8$ with
$y, \delta,d$ replaced by $y_1,\delta_1,d_1$.

For $j\leq y$, $\alpha,\beta\ll\frac{1}{\log y}$, \be
j^{\alpha},j^{\beta}=O(1),\ F(j,1+\alpha+\beta)=O(F_1(j,
1-2\delta)),\ee \be \frac{1}{F(j,1+\alpha)}=O(F_1(j, 1-2\delta)), \
\frac{1}{F(j,1+\beta)}=O(F_1(j, 1-2\delta)).\ee Then by \be
(1+\frac{1}{p^{1-2\delta}})^3\leq 1+\frac{7}{p^{1-2\delta}},\ee we
have \begin{eqnarray} U_2&=&O\Big(\frac{(\log\log
y)^6}{\log^4y}\sum_{j\leq
y}\frac{\mu^2(j)F_1^3(j,1-2\delta)}{j}\Big)\nonumber\\
&=&O\Big(\frac{(\log\log y)^6}{\log^4y}\sum_{j\leq
y}\frac{\mu^2(j)}{j}\sum_{n|j}\frac{d_7(n)}{n^{1-2\delta}}\Big)\nonumber\\
&=&O\Big(\frac{(\log\log y)^6}{\log^4y}\sum_{n\leq
y}\frac{d_7(n)}{n^{1-2\delta}}\sum_{j\leq y \atop
n|j}\frac{1}{j}\Big)\nonumber\\&=&O\Big(\frac{(\log\log
y)^6}{\log^4y}\sum_{n\leq
y}\frac{d_7(n)}{n^{2-2\delta}}\sum_{j_0\leq
y/n}\frac{1}{j_0}\Big)\nonumber\\&=&O\Big(\frac{(\log\log
y)^6}{\log^4y}\sum_{n=1}^{\infty}\frac{d_7(n)}{n^{2/3}}\log
y\Big)=O\Big(\frac{(\log\log y)^6}{\log^3y}\Big).\end{eqnarray}
Similarly,
\begin{eqnarray}
U_3&=&O\Big(\frac{(\log\log y)^4}{y^{2d}\log^2y}\sum_{j\leq
y}\frac{\mu^2(j)F_1^3(j,1-2\delta)}{j^{1-2d}}\Big)\nonumber\\
&=&O\Big(\frac{(\log\log y)^4}{y^{2d}\log^2y}\sum_{j\leq
y}\frac{\mu^2(j)}{j^{1-2d}}\sum_{n|j}\frac{d_7(n)}{n^{1-2\delta}}\Big)\nonumber\\
&=&O\Big(\frac{(\log\log y)^4}{y^{2d}\log^2y}\sum_{n\leq
y}\frac{d_7(n)}{n^{1-2\delta}}\sum_{j\leq y \atop
n|j}\frac{1}{j^{1-2d}}\Big)\nonumber\\&=&O\Big(\frac{(\log\log
y)^4}{y^{2d}\log^2y}\sum_{n\leq
y}\frac{d_7(n)}{n^{2-2d-2\delta}}\sum_{j_0\leq
y/n}\frac{1}{j_0^{1-2d}}\Big)\nonumber\\&=&O\Big(\frac{(\log\log
y)^4}{y^{2d}\log^2y}\sum_{n=1}^{\infty}\frac{d_7(n)}{n^{2/3}}(\frac{y^{2d}}{2d})\Big)=O\Big(\frac{(\log\log
y)^5}{\log^2y}\Big),\end{eqnarray}
\begin{eqnarray}
U_8&=&O\Big(\frac{(\log\log y)^5}{y^d\log^3y}\sum_{j\leq
y}\frac{\mu^2(j)F_1^3(j,1-2\delta)}{j^{1-d}}\Big)\nonumber\\
&=&O\Big(\frac{(\log\log
y)^5}{y^d\log^3y}\sum_{n=1}^{\infty}\frac{d_7(n)}{n^{2/3}}(\frac{y^d}{d})\Big)=O\Big(\frac{(\log\log
y)^6}{\log^3y}\Big).\end{eqnarray} For $j\leq y$, $0\leq m\leq I$,
it is easy to show \be G_m(\alpha,j)=O(\frac{1}{\log^{m+1}y}).\ee By
(4.10)-(4.12) and (4.16), we obtain \be U_4=O\Big(\frac{(\log\log
y)^3}{\log^3y}\sum_{j\leq
y}\frac{\mu^2(j)F_1^3(j,1-2\delta)}{j}\Big)=O\Big(\frac{(\log\log
y)^3}{\log^2y}\Big),\ee and similarly \be U_5=O\Big(\frac{(\log\log
y)^2}{y^d\log^2y}\sum_{j\leq
y}\frac{\mu^2(j)F_1^3(j,1-2\delta)}{j^{1-d}}\Big)=O\Big(\frac{(\log\log
y)^3}{\log^2y}\Big),\ee
 \be U_6=O\Big(\frac{(\log\log
y)^3}{\log^3y}\sum_{j\leq
y}\frac{\mu^2(j)F_1^3(j,1-2\delta)}{j}\Big)=O\Big(\frac{(\log\log
y)^3}{\log^2y}\Big),\ee \be U_7=O\Big(\frac{(\log\log
y)^2}{y^d\log^2y}\sum_{j\leq
y}\frac{\mu^2(j)F_1^3(j,1-2\delta)}{j^{1-d}}\Big)=O\Big(\frac{(\log\log
y)^3}{\log^2y}\Big).\ee We have the similar estimation for
$U_2',U_3',\cdots, U_8'$, then by (4.13)-(4.15), (4.17)-(4.20), and
note that $\frac{(\log\log y_1)^5}{\log^2y_1}=O(\frac{(\log\log
y)^5}{\log^2y})$, we obtain \be
U_2+\cdots+U_8+U_2'+\cdots+U_8'=O(\frac{(\log\log
y)^5}{\log^2y}).\ee

The rest of this section is due to evaluation of $U_1$. By Lemma 6
with $N=1$, Abel summation and (4.16), we have
\begin{eqnarray}
&&\sum_{j\leq
y}\frac{\mu^2(j)F(j,1+\alpha+\beta)}{jF(j,1+\alpha)F(j,1+\beta)}G_0(\alpha,j)G_0(\beta,j)\nonumber\\
&=&Y(\alpha,\beta)\int_1^y\frac{G_0(\alpha,\tau)G_0(\beta,\tau)}{\tau}d\tau+O(\frac{1}{\log^2
y}),\end{eqnarray} where \begin{eqnarray}
Y(\alpha,\beta)&=&\prod_p(1-\frac{1}{p^2})(1+\frac{1}{p+1}(\frac{1-p^{-1-\alpha-\beta}}{(1-p^{-1-\alpha})(1-p^{-1-\beta})}-1))\nonumber\\
&=&\prod_p(1-\frac{(p^{\alpha}-1)(p^{\beta}-1)}{(p^{1+\alpha}-1)(p^{1+\beta}-1)}).\end{eqnarray}
$Y(\alpha,\beta)$ is an analytic function in $|\alpha|<\frac14,\
|\beta|<\frac14$, then, \be
Y(\alpha,\beta)=Y(0,0)+O(\alpha)+O(\beta)=1+O(\frac{1}{\log y}).\ee
Let $t=\frac{\log\tau}{\log y_1}$, $\mu=\frac{\log x}{\log y_1}$ and
$\alpha=\frac{a}{\log T},\ \beta=\frac{b}{\log T},\ y=T^{\theta},
y_1=T^{\theta_1}$, then
\begin{eqnarray}
&&\int_1^y\frac{G_0(\alpha,\tau)G_0(\beta,\tau)}{\tau}d\tau\nonumber\\
&=&\log y_1\Big\{\int_0^1G_0(\alpha,e^{t\log y_1})G_0(\beta,e^{t\log
y_1})dt+\int_1^{\frac{\theta}{\theta_1}}G_0(\alpha,e^{t\log
y_1})G_0(\beta,e^{t\log
y_1})dt\Big\}\nonumber\\
&=&\log y_1\Big\{\int_0^1\Big(\alpha
P_1(1-\frac{\theta_1}{\theta}t)+\frac{1}{\log
y}P_1^{\prime}(1-\frac{\theta_1}{\theta}t)\nonumber\\&&\hspace*{1cm}+\sum_{l=2}^I\frac{(-1)^l}{(l-2)!\log
y_1}
\int_0^{1-t}P_l(1-t-\mu)\mu^{l-2}e^{-\mu\alpha\log y_1}d\mu\Big)\nonumber\\
&&\times\Big(\beta P_1(1-\frac{\theta_1}{\theta}t)+\frac{1}{\log
y}P_1^{\prime}(1-\frac{\theta_1}{\theta}t)\nonumber\\&&\hspace*{1cm}+\sum_{l=2}^I\frac{(-1)^l}{(l-2)!\log
y_1} \int_0^{1-t}P_l(1-t-\mu)\mu^{l-2}e^{-\mu\beta\log
y_1}d\mu\Big)dt\nonumber\\
&&+\int_1^{\frac{\theta}{\theta_1}}(\alpha
P_1(1-\frac{\theta_1}{\theta}t)+\frac{1}{\log
y}P_1^{\prime}(1-\frac{\theta_1}{\theta}t))(\beta
P_1(1-\frac{\theta_1}{\theta}t)+\frac{1}{\log
y}P_1^{\prime}(1-\frac{\theta_1}{\theta}t))dt\Big\}
\nonumber\\
&=&\frac{1}{\theta_1\log
T}\Big\{\int_0^1V_0(\theta,\theta_1,a,t)V_0(\theta,\theta_1,b,t)dt
\nonumber\\
&&\hspace*{2cm}+\int_1^{\frac{\theta}{\theta_1}}V_0^*(\theta,\theta_1,a,t)V_0^*(\theta,\theta_1,b,t)dt\Big\},\end{eqnarray}
where \begin{eqnarray}
 V_0(\theta,\theta_1,a,t)&=&a\theta_1
P_1(1-\frac{\theta_1}{\theta}t)+\frac{\theta_1}{\theta}P_1^{\prime}(1-\frac{\theta_1}{\theta}t)\nonumber\\
&&+\sum_{l=2}^I\frac{(-1)^l}{(l-2)!}
\int_0^{1-t}P_l(1-t-\mu)\mu^{l-2}e^{-a\theta_1\mu}d\mu,\end{eqnarray}
\be V_0^*(\theta,\theta_1,a,t)=a\theta_1
P_1(1-\frac{\theta_1}{\theta}t)+\frac{\theta_1}{\theta}P_1^{\prime}(1-\frac{\theta_1}{\theta}t).\ee
(4.22), (4.24), (4.25) and (4.16) gives
\begin{eqnarray}
&&\sum_{j\leq
y}\frac{\mu^2(j)F(j,1+\alpha+\beta)}{jF(j,1+\alpha)F(j,1+\beta)}G_0(\alpha,j)G_0(\beta,j)\nonumber\\
&=&\frac{1}{\theta_1\log
T}\Big\{\int_0^1V_0(\theta,\theta_1,a,t)V_0(\theta,\theta_1,b,t)dt
\nonumber\\
&&\hspace*{2cm}+\int_1^{\frac{\theta}{\theta_1}}V_0^*(\theta,\theta_1,a,t)V_0^*(\theta,\theta_1,b,t)dt\Big\}+O(\frac{1}{\log^2
y}).\end{eqnarray}

For $1\leq m\leq I$, let $j_0=j/p_1p_2\cdots p_m$, then
\begin{eqnarray}
&&\sum_{j\leq
y_1}\frac{\mu^2(j)F(j,1+\alpha+\beta)}{jF(j,1+\alpha)F(j,1+\beta)}G_0(\alpha,j)G_m(\beta,j)\sum_{p_1p_2\cdots
p_m|j}\log
p_1\log p_2\cdots\log p_m\nonumber\\
&=&\sum_{p_1p_2\cdots p_m\leq y_1}\frac{\mu^2(p_1p_2\cdots
p_m)F(p_1p_2\cdots p_m,1+\alpha+\beta)\log p_1\log p_2\cdots\log
p_m}{p_1p_2\cdots p_mF(p_1p_2\cdots p_m,1+\alpha)F(p_1p_2\cdots
p_m,1+\beta)}\times\nonumber
\\&&\sum_{j_0\leq
\frac{y_1}{p_1p_2\cdots p_m}\atop (j_0,p_1p_2\cdots
p_m)=1}\frac{\mu^2(j_0)F(j_0,1+\alpha+\beta)}{j_0F(j_0,1+\alpha)F(j_0,1+\beta)}G_0(\alpha,j_0p_1p_2\cdots
p_m)G_m(\beta,j_0p_1p_2\cdots p_m).\nonumber\\
&&\end{eqnarray} For $p_1p_2\cdots p_m\leq y_1$ and
$\mu(p_1p_2\cdots p_m)\not=0$, by lemma 6 with $N=p_1p_2\cdots p_m$,
Abel summation, (4.24) and (4.16), we obtain
\begin{eqnarray}
&&\sum_{j_0\leq \frac{y_1}{p_1p_2\cdots p_m}\atop (j_0,p_1p_2\cdots
p_m)=1}\frac{\mu^2(j_0)F(j_0,1+\alpha+\beta)}{j_0F(j_0,1+\alpha)F(j_0,1+\beta)}G_0(\alpha,j_0p_1p_2\cdots
p_m)G_m(\beta,j_0p_1p_2\cdots p_m)\nonumber\\
&&=Y(\alpha,\beta)\prod_{r=1}^m(1+\frac{1}{p_r})^{-1}(1+\frac{1}{p_r+1}(\frac{1-p_r^{-1-\alpha-\beta}}{(1-p_r^{-1-\alpha})(1-p_r^{-1-\beta})}-1))^{-1}\nonumber\\
&&\hspace*{1cm}\times \int_1^{\frac{y_1}{p_1p_2\cdots
p_m}}\frac{G_0(\alpha,\tau p_1p_2\cdots p_m)G_m(\beta,\tau
p_1p_2\cdots p_m)}{\tau}d\tau+O(\frac{\log\log y}{\log^{m+2}
y})\nonumber\\
&&=\prod_{r=1}^m\frac{(p_r^{1+\alpha}-1)(p_r^{1+\beta}-1)}{(p_r^{1+\alpha}-1)(p_r^{1+\beta}-1)+(p_r^{1+\alpha+\beta}-1)}\nonumber\\
&&\hspace*{1cm}\times \int_1^{\frac{y_1}{p_1p_2\cdots
p_m}}\frac{G_0(\alpha,\tau p_1p_2\cdots p_m)G_m(\beta,\tau
p_1p_2\cdots p_m)}{\tau}d\tau+O(\frac{\log\log y}{\log^{m+2}
y}).\end{eqnarray} Substitute (4.30) into (4.29), we get
\begin{eqnarray}
&&\hspace*{1cm}\sum_{j\leq
y_1}\frac{\mu^2(j)F(j,1+\alpha+\beta)}{jF(j,1+\alpha)F(j,1+\beta)}G_0(\alpha,j)G_m(\beta,j)\sum_{p_1p_2\cdots
p_m|j}\log
p_1\log p_2\cdots\log p_m\nonumber\\
&=&\sum_{p_1p_2\cdots p_m\leq y_1}\mu^2(p_1p_2\cdots
p_m)\prod_{r=1}^m\frac{(p_r^{1+\alpha+\beta}-1)\log
p_r}{(p_r^{1+\alpha}-1)(p_r^{1+\beta}-1)+(p_r^{1+\alpha+\beta}-1)}\nonumber
\\&&\hspace*{4cm}\times\int_1^{\frac{y_1}{p_1p_2\cdots
p_m}}\frac{G_0(\alpha,\tau p_1p_2\cdots p_m)G_m(\beta,\tau
p_1p_2\cdots p_m)}{\tau}d\tau\nonumber\\
&&+O(\frac{\log\log y}{\log^{m+2} y}\sum_{p_1p_2\cdots p_m\leq
y_1}\prod_{r=1}^m\frac{|p_r^{1+\alpha+\beta}-1|\log
p_r}{|(p_r^{1+\alpha}-1)(p_r^{1+\beta}-1)|})\nonumber\\&=&H_1+O(H_2),\end{eqnarray}
say. It is easy to show \be
\frac{p^{1+\alpha+\beta}-1}{(p^{1+\alpha}-1)(p^{1+\beta}-1)}=O(\frac{1}{p}),\ee
\be
\frac{p^{1+\alpha+\beta}-1}{(p^{1+\alpha}-1)(p^{1+\beta}-1)+(p^{1+\alpha+\beta}-1)}=\frac{1}{p}+O(\frac{1}{p^{3/2}})=O(\frac{1}{p}).\ee
Hence by lemma 4, \be H_2=O(\frac{\log\log y}{\log^{m+2}
y}\prod_{r=1}^m\sum_{p_r\leq y_1}\frac{\log
p_r}{p_r})=O(\frac{\log\log y}{\log^2y}),\ee and by (4.16) and lemma
5,\begin{eqnarray} H_1&=&\sum_{p_1p_2\cdots p_{m-1}\leq
y_1}\mu^2(p_1p_2\cdots
p_{m-1})\prod_{r=1}^{m-1}\frac{(p_r^{1+\alpha+\beta}-1)\log
p_r}{(p_r^{1+\alpha}-1)(p_r^{1+\beta}-1)+(p_r^{1+\alpha+\beta}-1)}\nonumber
\\&&\hspace*{5mm}\times\sum_{p_m\leq y_1/p_1p_2\cdots p_{m-1}\atop (p_m,p_1p_2\cdots p_{m-1})=1}\frac{\log p_m}{p_m}\int_1^{\frac{y_1}{p_1p_2\cdots
p_m}}\frac{G_0(\alpha,\tau p_1p_2\cdots p_m)G_m(\beta,\tau
p_1p_2\cdots p_m)}{\tau}d\tau\nonumber\\
&&\hspace*{5cm}+O(\prod_{r=1}^{m-1}\sum_{p_r\leq y_1}\frac{\log
p_r}{p_r}\sum_{p_m\leq y_1}\frac{\log
p_m}{p_m^{3/2}\log^{m+1}y})\nonumber\\&=&\sum_{p_1p_2\cdots
p_{m-1}\leq y_1}\mu^2(p_1p_2\cdots
p_{m-1})\prod_{r=1}^{m-1}\frac{(p_r^{1+\alpha+\beta}-1)\log
p_r}{(p_r^{1+\alpha}-1)(p_r^{1+\beta}-1)+(p_r^{1+\alpha+\beta}-1)}\nonumber
\\&&\hspace*{5mm}\times\sum_{p_m\leq y_1/p_1p_2\cdots p_{m-1}}\frac{\log p_m}{p_m}\int_1^{\frac{y_1}{p_1p_2\cdots
p_m}}\frac{G_0(\alpha,\tau p_1p_2\cdots p_m)G_m(\beta,\tau
p_1p_2\cdots p_m)}{\tau}d\tau\nonumber\\
&&+O(\sum_{p_1p_2\cdots p_{m-1}\leq y_1}\prod_{r=1}^{m-1}\frac{\log
p_r}{p_r}\sum_{p_m|p_1p_2\cdots p_{m-1}}\frac{\log
p_m}{p_m\log^{m+1}y})+O(\frac{1}{\log^2y})\nonumber\\&=&\sum_{p_1p_2\cdots
p_{m-1}\leq y_1}\mu^2(p_1p_2\cdots
p_{m-1})\prod_{r=1}^{m-1}\frac{(p_r^{1+\alpha+\beta}-1)\log
p_r}{(p_r^{1+\alpha}-1)(p_r^{1+\beta}-1)+(p_r^{1+\alpha+\beta}-1)}\nonumber
\\&&\hspace*{5mm}\times\sum_{p_m\leq y_1/p_1p_2\cdots p_{m-1}}\frac{\log p_m}{p_m}\int_1^{\frac{y_1}{p_1p_2\cdots
p_m}}\frac{G_0(\alpha,\tau p_1p_2\cdots p_m)G_m(\beta,\tau
p_1p_2\cdots p_m)}{\tau}d\tau\nonumber\\
&&+O(\frac{\log\log y}{\log^2y}).\end{eqnarray} We have by lemma 4,
Abel summation and (4.16)
\begin{eqnarray}
&&\sum_{p_m\leq y_1/p_1p_2\cdots p_{m-1}}\frac{\log
p_m}{p_m}\int_1^{\frac{y_1}{p_1p_2\cdots p_m}}\frac{G_0(\alpha,\tau
p_1p_2\cdots p_m)G_m(\beta,\tau
p_1p_2\cdots p_m)}{\tau}d\tau\nonumber\\
&=&\int_1^{\frac{y_1}{p_1p_2\cdots
p_{m-1}}}\frac{1}{\tau_m}\int_1^{\frac{y_1}{p_1p_2\cdots
p_{m-1}\tau_m}}\frac{G_0(\alpha,\tau p_1p_2\cdots
p_{m-1}\tau_m)G_m(\beta,\tau p_1p_2\cdots p_{m-1}\tau_m)}{\tau}d\tau
d\tau_m\nonumber\\
&&+O(\frac{1}{\log^{m+1}y}).
\end{eqnarray}
(4.31)-(4.36) and (4.16) gives
\begin{eqnarray}
&&\hspace*{1cm}\sum_{j\leq
y_1}\frac{\mu^2(j)F(j,1+\alpha+\beta)}{jF(j,1+\alpha)F(j,1+\beta)}G_0(\alpha,j)G_m(\beta,j)\sum_{p_1p_2\cdots
p_m|j}\log
p_1\log p_2\cdots\log p_m\nonumber\\
&=&\sum_{p_1p_2\cdots p_{m-1}\leq y_1}\mu^2(p_1p_2\cdots
p_{m-1})\prod_{r=1}^{m-1}\frac{(p_r^{1+\alpha+\beta}-1)\log
p_r}{(p_r^{1+\alpha}-1)(p_r^{1+\beta}-1)+(p_r^{1+\alpha+\beta}-1)}\times\nonumber
\\&&\int_1^{\frac{y_1}{p_1p_2\cdots
p_{m-1}}}\frac{1}{\tau_m}\int_1^{\frac{y_1}{p_1p_2\cdots
p_{m-1}\tau_m}}\frac{G_0(\alpha,\tau p_1p_2\cdots
p_{m-1}\tau_m)G_m(\beta,\tau p_1p_2\cdots p_{m-1}\tau_m)}{\tau}d\tau
d\tau_m\nonumber\\
&&+O(\frac{1}{\log^{m+1} y}\sum_{p_1p_2\cdots p_{m-1}\leq
y_1}\prod_{r=1}^{m-1}\frac{|p_r^{1+\alpha+\beta}-1|\log
p_r}{|(p_r^{1+\alpha}-1)(p_r^{1+\beta}-1)+(p_r^{1+\alpha+\beta}-1)|})+O(\frac{\log\log
y}{\log^2y})\nonumber\\&=&\sum_{p_1p_2\cdots p_{m-2}\leq
y_1}\mu^2(p_1p_2\cdots
p_{m-2})\prod_{r=1}^{m-2}\frac{(p_r^{1+\alpha+\beta}-1)\log
p_r}{(p_r^{1+\alpha}-1)(p_r^{1+\beta}-1)+(p_r^{1+\alpha+\beta}-1)}\nonumber\\
&&\hspace*{4cm}\times\sum_{p_{m-1}\leq y_1/p_1p_2\cdots
p_{m-2}}\frac{\log p_{m-1}}{p_{m-1}}\int_1^{\frac{y_1}{p_1p_2\cdots
p_{m-1}}}\frac{1}{\tau_m}\nonumber
\\&&\hspace*{1cm}\times\int_1^{\frac{y_1}{p_1p_2\cdots
p_{m-1}\tau_m}}\frac{G_0(\alpha,\tau p_1p_2\cdots
p_{m-1}\tau_m)G_m(\beta,\tau p_1p_2\cdots p_{m-1}\tau_m)}{\tau}d\tau
d\tau_m\nonumber\\
&&+O(\sum_{p_1p_2\cdots p_{m-2}\leq y_1}\prod_{r=1}^{m-2}\frac{\log
p_r}{p_r}\sum_{p_{m-1}|p_1p_2\cdots p_{m-2}}\frac{\log
p_{m-1}}{p_{m-1}\log^my})+O(\frac{\log\log
y}{\log^2y})\nonumber\\&=&\sum_{p_1p_2\cdots p_{m-2}\leq
y_1}\mu^2(p_1p_2\cdots
p_{m-2})\prod_{r=1}^{m-2}\frac{(p_r^{1+\alpha+\beta}-1)\log
p_r}{(p_r^{1+\alpha}-1)(p_r^{1+\beta}-1)+(p_r^{1+\alpha+\beta}-1)}\nonumber\\
&&\hspace*{4cm}\times\sum_{p_{m-1}\leq y_1/p_1p_2\cdots
p_{m-2}}\frac{\log p_{m-1}}{p_{m-1}}\int_1^{\frac{y_1}{p_1p_2\cdots
p_{m-1}}}\frac{1}{\tau_m}\nonumber
\\&&\hspace*{1cm}\times\int_1^{\frac{y_1}{p_1p_2\cdots
p_{m-1}\tau_m}}\frac{G_0(\alpha,\tau p_1p_2\cdots
p_{m-1}\tau_m)G_m(\beta,\tau p_1p_2\cdots p_{m-1}\tau_m)}{\tau}d\tau
d\tau_m\nonumber\\
&&+O(\frac{\log\log y}{\log^2y}).\end{eqnarray} Thus by induction,
lemma 7 and the same variable transformation as in (4.25), we get
\begin{eqnarray}
&&\hspace*{1cm}\sum_{j\leq
y_1}\frac{\mu^2(j)F(j,1+\alpha+\beta)}{jF(j,1+\alpha)F(j,1+\beta)}G_0(\alpha,j)G_m(\beta,j)\sum_{p_1p_2\cdots
p_m|j}\log
p_1\log p_2\cdots\log p_m\nonumber\\
&=&\int_1^{y_1}\frac{1}{\tau_1}\int_1^{\frac{y_1}{\tau_1}}\frac{1}{\tau_2}\cdots\int_1^{\frac{y_1}{\tau_1\tau_2\cdots
\tau_{m-1}}}\frac{1}{\tau_m}\nonumber\\&&\hspace*{1cm}\times\int_1^{\frac{y_1}{\tau_1\tau_2\cdots
\tau_m}}\frac{G_0(\alpha,\tau \tau_1\cdots \tau_m)G_m(\beta,\tau
\tau_1\cdots\tau_m)}{\tau}d\tau
d\tau_m\cdots d\tau_1+O(\frac{\log\log y}{\log^2y})\nonumber\\
&=&\int_1^{y_1}\frac{G_0(\alpha,\tau)G_m(\beta,\tau
)\log^m\tau}{m!\tau}d\tau+O(\frac{\log\log
y}{\log^2y})\nonumber\\
&=&\frac{1}{\theta_1\log
T}\int_0^1\frac{V_0(\theta,\theta_1,a,t)V_m(\theta,\theta_1,b,t)t^m}{m!}dt+O(\frac{\log\log
y}{\log^2y}),\end{eqnarray} and similarly
\begin{eqnarray}
&&\hspace*{1cm}\sum_{j\leq
y_1}\frac{\mu^2(j)F(j,1+\alpha+\beta)}{jF(j,1+\alpha)F(j,1+\beta)}G_m(\alpha,j)G_0(\beta,j)\sum_{p_1p_2\cdots
p_m|j}\log
p_1\log p_2\cdots\log p_m\nonumber\\
&=&\frac{1}{\theta_1\log
T}\int_0^1\frac{V_m(\theta,\theta_1,a,t)V_0(\theta,\theta_1,b,t)t^m}{m!}dt+O(\frac{\log\log
y}{\log^2y}),\end{eqnarray}
 with $V_0$ defined by
(4.26), \be V_1(\theta,\theta_1,a,t)=-2P_2(1-t)+\sum_{l=3}^I{\cal
C}_l^1\frac{(-1)^{l-1}}{(l-3)!}
\int_0^{1-t}P_l(1-t-\mu)\mu^{l-3}e^{-a\theta_1\mu}d\mu,\ee for
$2\leq m\leq I-2$,
\begin{eqnarray}
V_m(\theta,\theta_1,a,t)&=&a\theta_1
P_m(1-t)+P_m^{\prime}(1-t)-{\cal
C}_{m+1}^mP_{m+1}(1-t)\nonumber\\&+&\sum_{l=m+2}^I{\cal
C}_l^m\frac{(-1)^{l-m}}{(l-m-2)!}
\int_0^{1-t}P_l(1-t-\mu)\mu^{l-m-2}e^{-a\theta_1\mu}d\mu,\end{eqnarray}
and \be V_{I-1}(\theta,\theta_1,a,t)=a\theta_1
P_{I-1}(1-t)+P_{I-1}^{\prime}(1-t)-IP_I(1-t),\ee \be
V_I(\theta,\theta_1,a,t)=a\theta_1 P_I(1-t)+P_I^{\prime}(1-t).\ee

For $1\leq m_1,\ m_2\leq I$, by Lemma 8 we have
\begin{eqnarray}
&&\sum_{j\leq
y_1}\frac{\mu^2(j)F(j,1+\alpha+\beta)}{jF(j,1+\alpha)F(j,1+\beta)}G_{m_1}(\alpha,j)\sum_{p_1p_2\cdots
p_{m_1}|j}\log p_1\log p_2\cdots\log
p_{m_1}\nonumber\\&&\hspace*{5cm}\times
G_{m_2}(\beta,j)\sum_{p_1p_2\cdots p_{m_2}|j}\log
p_1\log p_2\cdots\log p_{m_2}\nonumber\\
&=&\sum_{k=0}^{\min(m_1,m_2)}{\cal P}_{m_1}^k{\cal
C}_{m_2}^k\sum_{j\leq
y_1}\frac{\mu^2(j)F(j,1+\alpha+\beta)}{jF(j,1+\alpha)F(j,1+\beta)}G_{m_1}(\alpha,j)G_{m_2}(\beta,j)\nonumber\\&&\hspace*{1cm}\times\sum_{p_1p_2\cdots
p_{m_1+m_2-k}|j}\log^2p_1\log^2p_2\cdots\log^2p_k\log
p_{k+1}\cdots\log p_{m_1+m_2-k}.\end{eqnarray} Similar to the proof
of (4.38), we have
\begin{eqnarray}
&&\sum_{j\leq
y_1}\frac{\mu^2(j)F(j,1+\alpha+\beta)}{jF(j,1+\alpha)F(j,1+\beta)}G_{m_1}(\alpha,j)G_{m_2}(\beta,j)\nonumber\\&&\hspace*{1cm}\times\sum_{p_1p_2\cdots
p_{m_1+m_2-k}|j}\log^2p_1\log^2p_2\cdots\log^2p_k\log
p_{k+1}\cdots\log p_{m_1+m_2-k}\nonumber\\
&=&\int_1^{y_1}\frac{\log\tau_1}{\tau_1}\cdots\int_1^{\frac{y_1}
{\tau_1\tau_2\cdots
\tau_{k-1}}}\frac{\log\tau_k}{\tau_k}\int_1^{\frac{y_1}{\tau_1\tau_2\cdots
\tau_{k}}}\frac{1}{\tau_{k+1}}\cdots\nonumber\\&&\int_1^{\frac{y_1}
{\tau_1\tau_2\cdots
\tau_{m_1+m_2-k}}}\frac{G_{m_1}(\alpha,\tau \tau_1\cdots
\tau_{m_1+m_2-k})G_{m_2}(\beta,\tau
\tau_1\cdots\tau_{m_1+m_2-k})}{\tau}d\tau
d\tau_{m_1+m_2-k}\cdots d\tau_1\nonumber\\
&&+O(\frac{\log\log y}{\log^2y}).\end{eqnarray}  Thus by lemma 9 and
the same variable transformation as in (4.25), we obtain
\begin{eqnarray}
&&\sum_{j\leq
y_1}\frac{\mu^2(j)F(j,1+\alpha+\beta)}{jF(j,1+\alpha)F(j,1+\beta)}G_{m_1}(\alpha,j)G_{m_2}(\beta,j)\nonumber\\&&\hspace*{1cm}\times\sum_{p_1p_2\cdots
p_{m_1+m_2-k}|j}\log^2p_1\log^2p_2\cdots\log^2p_k\log
p_{k+1}\cdots\log
p_{m_1+m_2-k}\nonumber\\&=&\int_1^{y_1}\frac{G_{m_1}(\alpha,\tau)G_{m_2}(\beta,\tau
)\log^{m_1+m_2}\tau}{(m_1+m_2)!\tau}d\tau+O(\frac{\log\log
y}{\log^2y})\nonumber\\
&=&\frac{1}{\theta_1\log
T}\int_0^1\frac{V_{m_1}(\theta,\theta_1,a,t)V_{m_2}(\theta,\theta_1,b,t)t^{m_1+m_2}}{(m_1+m_2)!}dt+O(\frac{\log\log
y}{\log^2y}).\end{eqnarray} with $V_1,\cdots V_I$ defined by
(4.40)-(4.43).

Substitute (4.28), (4.38), (4.39), (4.44) and (4.46) into (4.2), we
get \be U_1=\frac{1}{\theta_1\log T}\Big\{\int_0^1{\cal
F}(\theta,\theta_1,a, b, t)dt+\int_1^{\frac{\theta}{\theta_1}}{\cal
F}^*(\theta,\theta_1,a, b, t)dt\Big\}+O(\frac{\log\log
y}{\log^2y}),\ee where \be {\cal F}(\theta,\theta_1,a, b,
t)=\sum_{m_1=0}^I\sum_{m_2=0}^I\sum_{k=0}^{\min(m_1,m_2)}{\cal
P}_{m_1}^k{\cal
C}_{m_2}^k\frac{V_{m_1}(\theta,\theta_1,a,t)V_{m_2}(\theta,\theta_1,b,t)t^{m_1+m_2}}{(m_1+m_2)!},\ee
\be {\cal F}^*(\theta,\theta_1,a, b,
t)=V_0^*(\theta,\theta_1,a,t)V_0^*(\theta,\theta_1,b,t).\ee

(4.47) with (4.1), (4.21) gives \be
\Sigma(\alpha,\beta)=\frac{1}{\theta_1\log T}\Big\{\int_0^1{\cal
F}(\theta,\theta_1,a, b, t)dt+\int_1^{\frac{\theta}{\theta_1}}{\cal
F}^*(\theta,\theta_1,a, b, t)dt\Big\}+O(\frac{(\log\log
y)^5}{\log^2y}).\ee

 \vskip 8mm
\def\theequation{5.\arabic{equation}}
\setcounter{equation}{0} \centerline{\large\bf 5.\ Proof of Theorem
1} \vskip 5mm

By (4.50) and lemma 2, we have as $T\rightarrow\infty$,
\begin{eqnarray} g(\alpha,\beta,w)&=&\int_0^1\frac{{\cal
F}(\theta,\theta_1,b, a, t)-e^{-a-b}{\cal F}(\theta,\theta_1,-a, -b,
t)}{\theta_1(a+b)}dt\nonumber\\
&&+\int_1^{\frac{\theta}{\theta_1}}\frac{{\cal
F}^*(\theta,\theta_1,b, a, t)-e^{-a-b}{\cal F}^*(\theta,\theta_1,-a,
-b, t)}{\theta_1(a+b)}dt+o_{\delta}(1) \end{eqnarray} uniformly for
$a,\ b\ll 1$ and $T\leq w\leq 2T$. Then
\begin{eqnarray}
&&\frac{1}{\Delta\pi^{\frac12}}\int_{-\infty}^{\infty}e^{-(t-w)^2\Delta^{-2}}|BV(\sigma_0+it)|^2dt\nonumber\\
&=&Q(\frac{-d}{da})Q(\frac{-d}{db})g(\alpha,\beta,w)\Big|_{a=b=-R}\nonumber\\
&=&\int_0^1Q(\frac{-d}{da})Q(\frac{-d}{db})\Big(\frac{{\cal
F}(\theta,\theta_1,b, a, t)-e^{-a-b}{\cal F}(\theta,\theta_1,-a, -b,
t)}{\theta_1(a+b)}\Big)\Big|_{a=b=-R}dt\nonumber\\
&&+\int_1^{\frac{\theta}{\theta_1}}Q(\frac{-d}{da})Q(\frac{-d}{db})\Big(\frac{{\cal
F}^*(\theta,\theta_1,b, a, t)-e^{-a-b}{\cal F}^*(\theta,\theta_1,-a,
-b,
t)}{\theta_1(a+b)}\Big)\Big|_{a=b=-R}dt\nonumber\\
&&+o_{\delta}(1)\end{eqnarray} uniformly for $T\leq w\leq 2T$, with
$\Delta=T^{1-\delta}$. Thus it follows exactly as in section 3 of
Balasubramanian, Conrey and Heath-Brown [2] that
\begin{eqnarray} &&\frac{1}{T}\int_1^T|BV(\sigma_0+it)|^2dt\nonumber\\
&=&\int_0^1Q(\frac{-d}{da})Q(\frac{-d}{db})\Big(\frac{{\cal
F}(\theta,\theta_1,b, a, t)-e^{-a-b}{\cal F}(\theta,\theta_1,-a, -b,
t)}{\theta_1(a+b)}\Big)\Big|_{a=b=-R}dt\nonumber\\
&&+\int_1^{\frac{\theta}{\theta_1}}Q(\frac{-d}{da})Q(\frac{-d}{db})\Big(\frac{{\cal
F}^*(\theta,\theta_1,b, a, t)-e^{-a-b}{\cal F}^*(\theta,\theta_1,-a,
-b,
t)}{\theta_1(a+b)}\Big)\Big|_{a=b=-R}dt\nonumber\\
&&+o_{\delta}(1)\end{eqnarray}

Lemma 1 and (5.3) gives

\vskip 3mm

{\bf Theorem 2}. Let $T$ be a large parameter and $L=\log T$, $R$ be
a positive real number, $\theta<\frac47$, $\theta_1<\frac12$, $I\geq
2$ is a integer, $P_1$ is a real polynomial with $P_1(0)=0$ and
$P_1(1)=1$, $P_l (l=2,\cdots I)$ are real polynomials with
$P_l(0)=0$. Let $Q$ be a real polynomial with $Q(0)=1$ and
$Q'(x)=Q'(1-x)$. Then we have \begin{eqnarray} \kappa&&\geq
1-\frac{1}{R}\log(\int_0^1Q(\frac{-d}{da})Q(\frac{-d}{db})\Big(\frac{{\cal
F}(\theta,\theta_1,b, a, t)-e^{-a-b}{\cal F}(\theta,\theta_1,-a, -b,
t)}{\theta_1(a+b)}\Big)\Big|_{a=b=-R}dt\nonumber\\
&&+\int_1^{\frac{\theta}{\theta_1}}Q(\frac{-d}{da})Q(\frac{-d}{db})\Big(\frac{{\cal
F}^*(\theta,\theta_1,b, a, t)-e^{-a-b}{\cal F}^*(\theta,\theta_1,-a,
-b, t)}{\theta_1(a+b)}\Big)\Big|_{a=b=-R}dt), \end{eqnarray} with
${\cal F}(\theta,\theta_1,a, b, t)$ and ${\cal
F}^*(\theta,\theta_1,a, b, t)$ defined by (4.48) and (4.49),
$V_0,V_0^*,V_1,\cdots V_I$ defined by (4.26), (4.27), (4.40)-(4.43).

\vskip 3mm

{\bf Proof of theorem 1}. By theorem 2, with
\begin{eqnarray}
\theta&=&\frac47-\varepsilon,\ \ \theta_1=\frac12-\varepsilon,\ \ R=1.3025, \ \ I=3,\nonumber\\
P_1(x)&=&x+0.2950x(1-x)-2.2345x(1-x)^2+1.882x(1-x)^3,\nonumber\\
P_2(x)&=&0.0849x+1.9824x^2,\nonumber\\
P_3(x)&=&0.7516x,\nonumber\\
Q(x)&=&1-0.6684x-1.0798(\frac{x^2}{2}-\frac{x^3}{3})-5.0447(\frac{x^3}{3}-\frac{x^4}{2}+\frac{x^5}{5}),\nonumber
\end{eqnarray}
and let $\varepsilon\rightarrow 0$, we get
$$ \kappa\geq 0.4128.$$

\vskip 3mm

{\bf Corollary 1}. \ The '$\theta=1$' conjecture implies \be
\kappa\geq 0.6107.\ee

Here the '$\theta=1$' conjecture means that lemma 2, and then
theorem 2, is valid for any $\theta=\theta_1<1$.

\vskip 3mm

{\bf Proof}. \ With
\begin{eqnarray}
\theta&=&\theta_1=1-\varepsilon,\ \ R=0.7721, \ \ I=3,\nonumber\\
P_1(x)&=&x+0.1560x(1-x)-1.4045x(1-x)^2-0.0662x(1-x)^3,\nonumber\\
P_2(x)&=&2.0409x+0.2661x^2,\nonumber\\
P_3(x)&=&-0.0734x,\nonumber\\
Q(x)&=&1-0.7721x-0.1901(\frac{x^2}{2}-\frac{x^3}{3})-3.9627(\frac{x^3}{3}-\frac{x^4}{2}+\frac{x^5}{5}),\nonumber
\end{eqnarray}
and let $\varepsilon\rightarrow 0$, we get
$$ \kappa\geq 0.6107.$$

\vskip 5mm

{\bf Acknowledgements.} The author would like to thank Professor
Conrey for helpful discussion and comments, especially for remind me
to verify that if $\theta_1$ is permitted to take
$\theta_1=\frac47-\varepsilon$ as in earlier version of the
manuscript.

\vskip 8mm \centerline {\large\bf References} \vskip 5mm

\re{[1]} R. J. Anderson, {\sl Simple zeros of the Riemann
zeta-function}, J. Number Theory. {\bf 17}(1983), 176-182.

\re{[2]} R. Balasubramanian, J. B. Conrey and D. R. Heath-Brown,
{\sl Asympototic mean square of the product of the Riemann
zeta-function and a Dirichlet polynomial}, J. reine angew. Math.
{\bf 357}(1985), 161-181.

\re{[3]} H. M. Bui, J. B. Conrey and M. P. Young, {\sl More that
$41\%$ of the zeros of the zeta function are on the critical line},
http://arxiv.org/abs/1002.4127v2.

\re{[4]} J. B. Conrey, {\sl Zeros of derivatives of the Riemann's
$\xi$-function on the critical line}, J. Number Theory. {\bf
16}(1983), 49-74.

\re{[5]} J. B. Conrey, {\sl Zeros of derivatives of the Riemann's
$\xi$-function on the critical line, II}, J. Number Theory. {\bf
17}(1983), 71-75.

\re{[6]} J. B. Conrey, {\sl More than two fifths of the zeros of the
Riemann zeta function are on the critical line}, J. reine angew.
Math. {\bf 399}(1989), 1-26.

\re{[7]} J. M. Deshouillers and H. Iwaniec, {\sl Kloosterman sums
and Fourier coefficients of cusp forms}, Invent. Math. {\bf
70}(1982), 219-288.

\re{[8]} J. M. Deshouillers and H. Iwaniec, {\sl Power mean values
of the Riemann zeta function II}, Acta Arith. {\bf 48}(1984),
305-312.

\re{[9]} D. W. Farmer, {\sl Long mollifiers of the Riemann
zeta-function}, Mathematika. {\bf 40}(1993), 71-87.

\re{[10]} S. Feng, {\sl A note on the zeros of the derivative of the
 Riemann zeta function near the critical line}, Acta. Arith.
{\bf 120}(2005), 59-68.

\re{[11]} G. H. Hardy, {\sl Sur les z\'{e}ros de la fonction
$\zeta(s)$ de Riemann}, C. R. {\bf 158}(1914), 1012-1014.

\re{[12]} G. H. Hardy and J. E. Littlewood, {\sl The zeros of
Riemann's zeta-function on the critical line}, Math. Z. {\bf
10}(1921), 283-317.

\re{[13]} D. R. Heath-Brown, {\sl Simple zeros of the Riemann
zeta-function on the critical line}, Bull. London Math. Soc. {\bf
11}(1979), 17-18.

\re{[14]} N. Levinson, {\sl More than one third of zeros of
Riemann's zeta-function are on $\sigma=\frac12$}, Adv. Math. {\bf
13}(1974), 383-436.

\re{[15]} N. Levinson, {\sl Deduction of semi-optimal mollifier for
obtaining lower bounds for $N_0(T)$ for Riemann's zeta function},
Proc. Nat. Acad. Sci. USA. {\bf 72}(1975), 294-297.

\re{[16]} N. Levinson, {\sl Summing certain number theoretic sries
arising in the sieve}, J. Math. Anal. Appl. {\bf 22}(1968), 631-645.

\re{[17]} N. Levinson and H. L. Mongomery, {\sl zeros of the
derivative of the Riemann zeta-function}, Acta mathematica.{\bf
133}(1974), 49-65.

\re{[18]} S. Lou and Q. Yao, {\sl A lower bound for zeros of
Riemann's zeta function on the line $\sigma=\frac12$}, (chinese)
Acta Mathematica Sinica. {\bf 24}(1981), 390-400.

\re{[19]} H. L. Montgomery, {\sl Selberg's work on the
zeta-function}, Number theory, trace formulas and discrete groups,
Academic Press, Boston, (1989), 157-168.

\re{[20]} B. Riemann, {\sl \"{U}ber die Anzahl Primzahlen unter eine
gegebener Gr\"{o}sse}. Monatsber. Akad. Berlin. (1859), 671-680.

\re{[21]} A. Selberg, {\sl On the zeros of Riemann's zeta-function},
Skr. Norske Vid. Akad. Oslo. {\bf 10}(1942), 1-59.

\re{[22]} A. Selberg, {\sl On an elementary method in the theory of
primes}, Norske Vid. Selsk. Forh. {\bf 19}(1947), 64-67.

\re{[23]} C. L. Siegel, {\sl \"{U}ber Riemann's Nachlass zur
analytischen Zahlentheorie}, Quellen und Studien zur Geschichte der
Math. Astr. und Physik, Abt. B: Studien {\bf 2}(1932), 45-80.

\re{[24]} E. C. Titchmarsh, {\sl The Theory of the Riemann Zeta
Function}. 2nd edition. revised by D. R. Heath-Brown, Clarendon
Press, Oxford, 1986.

\re{[25]} Y. Zhang, {\sl On the zeros of $\zeta^{\prime}(s)$ near
the critical line}, Duke Math. J. {\bf 110}(2001), 555-571.

\end{document}